\documentclass[11pt]{article}
\usepackage{amsmath,amssymb,amsthm,amsfonts}
\usepackage{mathtools}
\usepackage{graphicx}
\usepackage[round,longnamesfirst]{natbib}
\usepackage{hyperref}
\usepackage{multimedia}
\usepackage{bm}
\usepackage{proba}
\usepackage{enumerate}
\usepackage{natbib}
\usepackage{lscape}
\usepackage{placeins}
\usepackage{longtable}
\usepackage{epstopdf}
\usepackage{graphics}
\usepackage{float}
\usepackage{placeins}
\usepackage{longtable}
\usepackage{appendix}
\usepackage{caption}
\usepackage{subcaption}
\usepackage{bigints}
\usepackage{booktabs}
\usepackage{xr}
\usepackage{xcolor}
\usepackage{geometry}
\usepackage{changepage}

\renewcommand{\cite}{\citet}
 \bibpunct{(}{)}{,}{a}{,}{,}
\newtheorem{assum}{Assumption}[section]

\newtheorem{thm}{Theorem}[section]
\newtheorem{rem}{Remark}[section]
\definecolor{dgreen}{rgb}{0,0.5,0}
\definecolor{dblue}{rgb}{0,0,0.9}
\definecolor{dred}{rgb}{0.6,0.0,0.1}
\definecolor{dgold}{rgb}{0.5,0.3,0.0}
\definecolor{dvio}{rgb}{0.6,0.3,0.5}
\definecolor{gray}{rgb}{0.5,0.5,0.5}

\newcommand{\bee}{\begin{equation}}
\newcommand{\eee}{\end{equation}}

\newcommand{\wh}{\widehat}

\newcommand{\wtl}{\widetilde}

\newcommand{\EE}{\mathbf{E}}

\newcommand{\arginf}{\mathrm{arginf}}
\newcommand{\argmax}{\mathrm{argmax}}
\renewcommand{\E}{\mathbf{E}}
\makeatletter \@addtoreset{equation}{section} \makeatother \renewcommand{\theequation}{\thesection.\arabic{equation}}
\bibpunct{(}{)}{,}{a}{,}{,}
\allowdisplaybreaks
\newtheorem{example1}{Example}
\newenvironment{contexample1}{
   \addtocounter{example1}{-1} \begin{example1}[continued]}{
   \end{example1}}
\newtheorem{example2}{Example}

\newtheorem{example3}{Example}

\newtheorem{example4}{Example}
\newenvironment{contexample4}{
   \addtocounter{example4}{-1} \begin{example4}[continued]}{
   \end{example4}}
\newtheorem{example5}{Example}

\newtheorem{df}{Definition}[section]
\input{tcilatex}
\begin{document}

\title{\vspace*{-.1in}\textbf{Bayesian Estimation and Comparison of \\
Conditional Moment Models}\footnote{We are enormously grateful to the editor,
associate editor and referee for their constructive comments and penetrating questions that spawned many improvements. Disclaimer: The views expressed in this paper are solely those of the authors and do not necessarily reflect the views of the Federal Reserve Bank of Philadelphia or the Federal Reserve System.}  }
\author{Siddhartha Chib\thanks{%
Olin Business School, Washington University in St. Louis, Campus Box 1133, 1
Brookings Drive, St. Louis, MO 63130. e-mail: chib@wustl.edu.} \quad \quad
Minchul Shin\thanks{%
Federal Reserve Bank of Philadelphia, Ten Independence Mall, Philadelphia,
PA 19106. e-mail: visiblehand@gmail.com. } \quad \quad Anna Simoni\thanks{%
CREST, CNRS, \'{E}cole Polytechnique, ENSAE, 5, Avenue Henry Le Chatelier,
91120 Palaiseau - France, e-mail: simoni.anna@gmail.com.} }
\date{December 2019, December 2020, July 2021}
\maketitle

\begin{abstract}
We consider the Bayesian analysis of models in which the unknown
distribution of the outcomes is specified up to a set of conditional moment
restrictions. The nonparametric exponentially tilted empirical likelihood
function is constructed to satisfy a sequence of unconditional
moments based on an increasing (in sample size) vector of approximating
functions (such as tensor splines based on the splines of each conditioning
variable). For any given sample size, results are robust to the number of
expanded moments. We derive Bernstein-von Mises theorems for
the behavior of the posterior distribution under both
correct and incorrect specification of the conditional moments,
subject to growth rate conditions (slower under misspecification)
on the number of approximating functions.
A large-sample theory for comparing different conditional moment
models is also developed. The central result is that the marginal likelihood
criterion selects the model that is less misspecified.
We also introduce sparsity-based model search for high-dimensional
conditioning variables, and provide efficient MCMC computations
for high-dimensional parameters. Along with clarifying examples, the
framework is illustrated with real-data applications to risk-factor
determination in finance, and causal inference under conditional ignorability.
\end{abstract}



\thispagestyle{empty}

\setlength{\baselineskip}{21.0pt}



\section{Introduction}

\label{Intro}

\sloppy 

We tackle the problem of prior-posterior inference when the only available
information about the unknown parameter $\theta \in \Theta \subset \mathbb{R}^{p}$ is supplied by a set of \textit{conditional moment} (CM) restrictions
\begin{equation}
\mathbf{E}^{P}[\rho (X,\theta )|Z]=0,  \label{model_introduction}
\end{equation}%
where $\rho (X,\theta )$ is a $d$-vector of known functions of a $\mathbb{R}^{d_{x}}$-valued random vector $X$ and the unknown $\theta $, and $P$ is the
unknown conditional distribution of $X$ given a $\mathbb{R}^{d_{z}}$-valued
random vector $Z$. Such models are important because many standard models in
statistics can be recast in terms of CM restrictions. These models also arise
naturally in causal inference, missing data problems, and in models derived
from theory in economics and finance. Because the CM conditions
constrain the set of possible distributions $P$, we say that the model is
correctly specified if the true data generating process $P_{\ast }$
is in the set of distributions constrained to satisfy these moment conditions for some $\theta \in \Theta $, while the model is misspecified if $P_{\ast }$ is not in the set of
implied distributions for any $\theta \in \Theta $.

A different starting point is when one is given the \emph{unconditional} moments,
say $\mathbf{E}^{P}[g(X,\theta )]=0$. Prior-posterior analysis
can then be based on the empirical likelihood, for example, \citet*{Lazar} and many others, or the
exponentially tilted empirical likelihood (ETEL), as in \citet*{Schennach2005%
} and \citet*{ChibShinSimoni2018}. Developing a Bayesian framework for CM
models is important. While it is true that the conditional moments
imply that $\rho (X,\theta )$ is uncorrelated with $Z$, i.e., $\mathbf{E}^{P}[\rho (X, \theta ) \otimes
Z]=0$, where $\otimes $ is the Kronecker product operator,
the conditional moments assert even more, that $\rho (X,\theta )$ is uncorrelated with any measurable,
bounded function of $Z$. Thus, there is an efficiency loss if this information is ignored.

We approach this problem by first constructing  $K$ unconditional moments
\begin{equation}
\mathbf{E}^{P}[\rho (X,\theta )\otimes q^{K}(Z)]=0
\label{model_unconditional_introduction}
\end{equation}%
based on an increasing (in sample size) vector of approximating
functions, $q^{K}(Z)\coloneqq (q_{1}^{K}(Z),\ldots ,q_{K}^{K}(Z))'$, obtained, for instance, from splines of each variable in $Z$
\citep*{DonaldImbensNewey2003}. Efficiency loss is avoided as the number of moments increases
with sample size. Next, for each sample size and for each $\theta$, the nonparametric exponentially tilted empirical likelihood (ETEL) function is constructed to satisfy these unconditional moments.
Unlike the empirical likelihood, the ETEL function has a fully Bayesian interpretation. It is  the likelihood that emerges from integrating out $P$ with respect to a nonparametric prior that satisfies the CMs. The posterior of interest is then this nonparametric likelihood multiplied by a prior distribution of the parameters. Due to the fact that the nonparametric likelihood is limited to a set $H_{n,K}$ of $\theta$ values for which the empirical counterpart of the moment conditions \eqref{model_unconditional_introduction} are equal to 0, 
the posterior (equivalently, the prior) is truncated to the set $H_{n,K}$.

We study the prior-posterior mapping on many fronts, taking up the question of misspecified models, model comparisons, and computations, combining careful theoretical work with the needs of applications. The posterior distribution is shown to satisfy Bernstein-von Mises (BvM) theorems in both the correct and misspecified cases. In the former case the growth of $K$ (for approximating functions given by splines) is at most $n^{1/6}$, where $n$ is the sample size. The asymptotic posterior variance is then equal to the semiparametric efficiency bound derived in \citet*{Chamberlain1987}.
In the latter case, in parallel with \citet*{kleijn2012}, the posterior distribution of the centered and scaled parameter $\sqrt{n}(\theta -\theta _{\circ })$,
where $\theta _{\circ }$ is the pseudo-true value, converges to a Normal
distribution with variance that now is different from the variance of the
frequentist estimator. Interestingly, this convergence holds only if $K$ increases more slowly than in the correctly specified case. This can be interpreted as limiting the
number of implied unconditional moments to limit the magnification of the
misspecification.

We informally use these rate conditions from the theoretical analysis
to guide the range of choice of $K$ for any given $n$.
Due to the fact that for a fixed $n$ the volume (prior probability content)
of the region of truncation $H_{n,K}$ decreases with $K$ (a
result of more restrictions), values of $K$ beyond the range
recommended by the theory amplify the Bayesian bias,
and, hence, should be avoided. Large values of $K$ can also produce
rank-deficiency of the approximating functions basis matrix and,
in the event of a misspecified model, increase misspecification.
Around the values of $K$ we recommend, the posterior distribution is
generally robust to $K$, and little fine-tuning is necessary.

Finite sample summaries of the posterior distribution are
obtained by Markov chain Monte Carlo (MCMC) methods.
Since the posterior is underpinned by a non-parametric
likelihood, and the effective prior is truncated,
efficient sampling is not automatic. However,
after extensive study, we have produced a near-black-box MCMC
approach (available as a R-package) that is based on the
tailored Metropolis-Hastings (M-H) algorithm of \citet*{ChibGreenberg1995}
and its randomized version in \citet*{ChibRamamurthy2010}.

The entire paper is interspersed with examples of pedagogical
importance and practical relevance. Real data applications
to risk-factor determination in finance, and causal inference
under conditional ignorability, are included.

It is worth noting that previous Bayesian work on
conditional moments, for example,
\citet*{liao2011}, \citet*{FS2012JoE,
florens_simoni_2016}, \citet*{kato2013},
\citet*{chen2017monte} and \cite{LIAOSimoni2019}, has little overlap with the
discussion here. A major difference is that none of these
papers adopt the fully Bayesian ETEL framework. Another is that these papers
examine a different class of CM models. Finally, none of these papers
takes up the question of model comparisons. Nonetheless, these papers
and the current work, taken together, represent an important
broadening of the Bayesian enterprise to new classes of models.

The rest of the paper is organized as follows. Section \ref{Setting} has the
sketch of the conditional moment setting. Section \ref{s_3} discusses
the prior-posterior analysis and the large-sample properties of the
posterior distribution. Section \ref{s_4} is concerned with
the problem of comparing CM models via marginal likelihoods.
In Section \ref{s_Additional_Topics} two extensions are considered
and Section \ref{s_Applications} has real data
applications to finance and causal inference. Section \ref{s_6} concludes.
Proofs are in the online supplementary appendix.

\section{Setting and Motivation}

\label{Setting}

Let $X\coloneqq(X_{1}^{\prime },X_{z}^{\prime })^{\prime }$ be an $\mathbb{R}%
^{d_{x}}$-valued random vector and $Z\coloneqq(Z_{1}^{\prime },X_{z}^{\prime
})^{\prime }$ be an $\mathbb{R}^{d_{z}}$-valued random vector. The vectors $%
Z $ and $X$ have elements in common if the dimension of the subvector $X_{z}$
is non-zero. Moreover, we denote $W\coloneqq(X^{\prime },Z_{1}^{\prime })^{\prime
}\in \mathbb{R}^{d_{w}}$ and its (unknown) joint distribution by $P$. By
abuse of notation, let $P$ also denote the associated conditional
distribution. Suppose that we are given a random sample $%
W_{1:n}=(W_{1},\ldots ,W_{n})$ of $W$. Hereafter, $\mathbf{E}
^{P}[\cdot ]$ is the expectation with respect to $P$ and $\mathbf{E}%
^{P}[\cdot |\cdot ]$ is the conditional expectation with respect to the
conditional distribution associated with $P$.

The parameter of interest is $\theta \in \Theta \subset \mathbb{R}^{p}$,
which is related to the conditional distribution $P$ through the conditional
moment restrictions
\begin{equation}
\mathbf{E}^{P}[\rho (X,\theta )|Z]=0,  \label{eq:condmom}
\end{equation}%
where $\rho (X,\theta )$ is a $d$-vector of known functions. Many
interesting and important models in statistics fall into this framework.

\begin{example1}
(Linear model with heteroscedasticity of unknown form) Suppose that
\begin{equation}
\mathbf{E}^{P}[(Y-\theta _{0}-\theta _{1}X)|Z]=0,  \label{eq:ex01_condi1}
\end{equation}%
where $\rho (X,\theta )=\left( Y-\theta _{0}-\theta _{1}X\right) $, $Z=(1,X)$
and $d=1$. This CM model is consistent with the
data generating process (DGP) $Y=\theta _{0}+\theta _{1}X+\varepsilon $,
where $\varepsilon = h(X)U$, and $(X,U)$ (independent) follow some unknown
distribution $P$, with $E(U)=0$, and the heteroscedasticity
function $h(X)$ is unknown. The restrictions
\begin{equation}
\begin{split}
\mathbf{E}^{P}[\left( Y-\theta _{0}-\theta _{1}X\right) |Z] =0 \quad \text{and} \quad \mathbf{E}^{P}[(Y-\theta _{0}-\theta _{1}X)^{3}|Z] =0,
\end{split}
\label{eq_ex01_mom}
\end{equation}%
where now $\rho (X,\theta )$ is a $(2\times 1)$ vector of functions,
additionally impose that $\varepsilon $ is conditionally symmetric.
\end{example1}

Note that in the foregoing example, the two unconditional moment conditions
\begin{equation}
\mathbf{E}^{P}[(Y-\theta _{0}-\theta _{1}X)\otimes (1,X)^{\prime }]=0,
\label{eq:ex01 uncondmom}
\end{equation}%
which assert that: (i) $\varepsilon $ has mean zero and (ii) $\varepsilon$ is
uncorrelated with $X$, are weaker but, if the CM model is correct,
less informative about $\theta$.


\section{Prior-Posterior Analysis}

\label{s_3}


\subsection{Expanded Moment Conditions}

The starting point, as in the frequentist approaches of
\citet*{DonaldNewey2001}, \citet*{AiChen2003} and
\cite{carrasco_florens_2000}, is a transformation of the CM
restrictions into unconditional moment restrictions.
Following \citet*{DonaldImbensNewey2003}, let $q^{K}(Z) \coloneqq (q_{1}^{K}(Z),\ldots ,q_{K}^{K}(Z))^{\prime }$, $K>0$, denote
a $K$-vector of real-valued functions of $Z$, for instance, splines.
Suppose that these functions satisfy the following condition for the distribution $P$.

\begin{assum}
\label{Ass_1_DIN} For all $K$, $\mathbf{E}^P[q^K(Z)^{\prime} q^K(Z)]$ is
finite, and for any function $a(z): \R^{d_{z}} \rightarrow \R $ with $%
\mathbf{E}^P[a(Z)^2]< \infty$ there are $K\times 1$ vectors $\gamma_K$ such
that as $K\rightarrow \infty$,
\begin{equation*}
\mathbf{E}^P[(a(Z)-q^K(Z)^{\prime }\gamma_K)^2]\rightarrow 0.
\end{equation*}
\end{assum}
\noindent Now, let $\theta_*$ be the value of $\theta$ that satisfies \eqref{eq:condmom} for the true $P$. If $\mathbf{E}^{P}[\rho (X,\theta_*)^{\prime_*}\rho (X,\theta )]<\infty $,
then \citet*[Lemma 2.1]{DonaldImbensNewey2003} established that: (1) if
equation \eqref{eq:condmom} is satisfied with $\theta =\theta_{\ast }$, then
$\mathbf{E}^{P}[\rho (X,\theta _{\ast })\otimes q^{K}(Z)]=0$ for all $K$;
(2) if equation \eqref{eq:condmom} is not satisfied by $\theta =\theta_{\ast }$, then $\mathbf{E}%
^{P}[\rho (X,\theta _{\ast })\otimes q^{K}(Z)]\neq 0$, for all large enough $%
K$.

Henceforth, we let
$g(W,\theta )\coloneqq\rho (X,\theta )\otimes q^{K}(Z)$
denote the \textit{expanded functions} and refer to
\begin{equation}
\mathbf{E}^{P}[g(W,\theta )]=0,  \label{eq:uncondmom}
\end{equation}
as the \textit{expanded moments}. Under the stated assumptions, the expanded moments
are equivalent to the CM restrictions \eqref{eq:condmom}, as $K\rightarrow \infty$.

In our numerical examples, we construct $q^{K}(Z)$
using the natural cubic spline basis of \citet*{Chib2010},
with $K$ fixed at a given value, as in sieve estimation.
If $Z$ consists of more than one element, say
$(Z_{1},Z_{2},Z_{3})$ where $Z_{1}$ and $Z_{2}$ are continuous variables and
$Z_{3}$ is binary, then the basis matrix $B$ is constructed as follows.
Let $\boldsymbol{z}_{j}$ denote the $n\times 1$ sample data on $Z_{j}$ $\left(
j\leq 3\right) $. Let $\boldsymbol{Z}=\left( \boldsymbol{z}_{1},\boldsymbol{z%
}_{2},\boldsymbol{z}_{1}\odot \boldsymbol{z}_{2},\boldsymbol{z}_{1}\odot
\boldsymbol{z}_{3},\boldsymbol{z}_{2}\odot \boldsymbol{z}_{3}\right) $
denote the $n\times 5$ matrix of the continuous data and interactions of the
continuous data and the binary data. Now suppose $\left( \tau _{j1},\ldots
,\tau _{jK}\right) $, for $j=1,\ldots,5$ are $K$ knots based on each column of $\boldsymbol{Z}$
and let $B_{j}$ denote the corresponding $n\times K$ matrix of cubic spline
basis functions. Then, $B$ is given by
\begin{equation*}
B=\left[ B_{1}\;\vdots \;B_{2}^{\ast }\;\vdots \;B_{3}^{\ast }\;\vdots
\;B_{4}^{\ast }\;\vdots \;B_{5}^{\ast }\;\vdots \;\boldsymbol{Z}_{3}\right] ,
\end{equation*}%
where $B_{j}^{\ast }$ $\left( j=2,3,4,5\right) $ is the $n\times (K-1)$
matrix in which each column of $B_{j}$ is subtracted from its first and then
the first column is dropped, see \citet*{Chib2010}. Thus, the dimension of
this $B$ matrix is $n\times K^{\ast}$, where $K^{\ast} = \left( 5K-4+1\right)$. If
$K^{\ast}$ is large, in relation to $n$, data-compression methods can be employed.
Specifically, let $R$ denote the $K^{\ast} \times K^{\ast}$ orthogonal matrix of
eigenvectors from the singular value decomposition of $B$, and let
$e$ denote the corresponding $K^{\ast} \times 1$ vector of eigenvalues.
Then, after employing the rotation $BR$, the columns of $BR$ corresponding
to small values of $e$ are dropped, and the resulting column-reduced $BR$
matrix is taken as the basis matrix. We refer to this as the \textit{rotated
column reduced} basis matrix. To define the expanded
functions, let $\rho _{l}(\boldsymbol{X},\theta )$ $(l\leq d)$ denote a $%
n\times 1$ vector of the $l$th element of $\rho (X,\theta )$ evaluated at
the sample data matrix $\boldsymbol{X}$. Then, the expanded functions for
the sample observations are obtained by multiplying $\rho _{l}(\boldsymbol{X}%
,\theta )$ by the matrix $B$ (or by the rotated column reduced $B$)
and concatenating. We use versions of this approach in our examples.

\subsection{Posterior distribution}

We base the prior-posterior mapping, for each sample size and $\theta$,
on the nonparametric exponentially tilted empirical likelihood (ETEL) function.
The ETEL has a fully Bayesian interpretation \citep*{Schennach2005}
as an integrated likelihood, integrated over the prior on the data distribution
$P$ that satisfies the given moments. Other such priors exist,
for example, \citet*{KitamuraOtsu2011}, \citet*{Shin2014}
and \citet*{FlorensSimoni2015}, that lead to different integrated likelihoods.

The ETEL function takes the form
\begin{equation}
p(W_{1:n}|\theta ,K)=\prod_{i=1}^{n}\widehat{p}_{i}(\theta ),
\end{equation}%
where $\{\widehat{p}_{i}(\theta ),i=1,\ldots
,n\}$ are the probabilities that minimize the Kullback-Leibler divergence between the
probabilities $(p_{1},\ldots ,p_{n})$ assigned to each sample observation
and the empirical probabilities $(\frac{1}{n},\ldots ,\frac{1}{n})$, subject
to the conditions that the probabilities $(p_{1},\ldots ,p_{n})$ sum to one
and that the expectation under these probabilities satisfy the given
unconditional moment conditions \eqref{eq:uncondmom}.

Specifically, $\{\widehat{p}%
_{i}(\theta ),i=1,\ldots ,n\}$ are the solution of the following problem:
\begin{equation}
\max_{p_{1},\ldots ,p_{n}}\sum_{i=1}^{n}\left[ -p_{i}\log (np_{i})\right]
\qquad \text{subject to: }\sum_{i=1}^{n}p_{i}=1,\quad
\sum_{i=1}^{n}p_{i}g(w_{i},\theta )=0,\quad p_{i}\geq 0
\end{equation}%
(see \citet*{Schennach2005} for a proof). In practice, the solution of this problem
emerges from the dual (saddlepoint)
representation (see \textit{e.g.} \citet*{csiszar1984}) as
\begin{equation}
\widehat{p}_{i}(\theta )\coloneqq\frac{e^{\widehat{\lambda }(\theta )^{\prime
}g(w_{i},\theta )}}{\sum_{j=1}^{n}e^{\widehat{\lambda }(\theta )^{\prime
}g(w_{j},\theta )}}, \qquad i=1,\ldots ,n,  \label{eq_ETEL_weights}
\end{equation}%
where $\widehat{\lambda }(\theta )=\arg \min_{\lambda \in \mathbb{R}^{dK}}%
\frac{1}{n}\sum_{i=1}^{n}e^{\lambda ^{\prime }g(w_{i},\theta )}$ is the
estimated tilting parameter.

Let $\mathcal{C}o(\theta) \coloneqq \left\{ \sum_{i=1}^{n}p_{i}g(w_{i},\theta),p_{i}\geq 0, \sum_{i=1}^n p_i = 1\right\}$ be the convex hull of $\{g(w_{i},\theta )\}_{i=1}^n$ for a given $\theta$ and $\overline{\mathcal{C}o(\theta)}$ denote its interior. Let $H_{n,K}\coloneqq\left\{ \theta \in \Theta ; 0\in \overline{\mathcal{C}o(\theta)}\right\} $ denote the set
of $\theta $ values for which the empirical moment conditions hold. Then, the posterior
distribution is the truncated distribution given by
\begin{equation}
\pi (\theta |w_{1:n},K)\propto \pi (\theta )\prod_{i=1}^{n}\frac{e^{\widehat{%
\lambda }(\theta )^{\prime }g(w_{i},\theta )}}{\sum_{j=1}^{n}e^{\widehat{%
\lambda }(\theta )^{\prime }g(w_{j},\theta )}}1\{\theta\in H_{n,K}\},
\label{eq_betel posterior}
\end{equation}%
where $1\{\cdot\}$ is the indicator function.

Combining the indicator function with the prior, we see that
the (effective) prior is truncated to $\theta\in H_{n,K}$. This fact can be used to argue
that, for fixed $n$, it is not desirable to have a large $K$. This is because as $K$
increases for a given $n$, the support of the prior shrinks (equivalently, the prior probability
content of the region of truncation decreases), due to the fact that more restrictions
are imposed. We refer to this prior probability content by the shorthand, volume. Reduction in the volume tends to increase the Bayesian bias and reduce the posterior spread, without any change in the data, with deleterious impact on the posterior. In practice, we use the rule $2n^{1/6}$ to fix $K$. Larger values than this can, of course, be tried, but one should make sure that the volume of
$H_{n,K}$ does not become much smaller than one. Around the values of $K$ we recommend, the posterior distribution is generally robust to $K$, and little fine-tuning is necessary.

\begin{contexample1}
To illustrate the role of $K$ in the prior-posterior analysis, and its impact on the volume (prior probability content) of $H_{n,K}$, we create a set of simulated data $\{y_{i},x_{i}\}_{i=1}^{n}$, $n=250$, with covariates $X\sim \mathcal{U}(-1,2.5)$, intercept $\theta _{0}=1$, slope $\theta _{1}=1$, and $\varepsilon _{i}$ is distributed according to $\varepsilon _{i}\sim\mathcal{SN}(m(x_{i}),h(x_{i}),s(x_{i}))$, where $\mathcal{SN}(m,h,s)$ is
the skew normal distribution with location, scale, and shape parameters
given by $(m,h,s)$, each depending on $x_{i}$. When $s$ is zero, $%
\varepsilon _{i}$ is normal with mean $m$ and standard deviation $h$. We set
$m(x_{i})=-h(x_{i})\sqrt{2/\pi }s(x_{i})/(\sqrt{1+s(x_{i})^{2}})$,
so that $\EE^{P}[\varepsilon |X]=0$.

Suppose that $h(x)=\sqrt{\exp (1+0.7x+0.2x^{2})}$ and $s(x)=1+x^{2}$.
Model parameters are estimated solely from the condition
$\EE^{P}[\varepsilon |Z]=0$, $Z=(1,X)$. The prior is
the default independent student-$t$ distribution with location 0,
dispersion 5, and degrees of freedom 2.5, truncated to $H_{n,K}$.
The posterior is computed for $K$ given by 2, $2n^{1/6}$, 9 and 20
(the value $2n^{1/6}$ is based on the theory below for splines approximating functions).
The results are shown in Table \ref{Tab:1}. Importantly,
when $K$ is close to the value suggested by theory, the posterior
distribution is robust to $K$. However,
when $K=20$, quite different from the recommended value,
the Bayesian bias is larger and the posterior
standard deviation is smaller, without any change in the data. This is
due to the effect of the prior, in the following way.
As $K$ increases for a fixed $n$, the volume of $H_{n,K}$ decreases,
equivalently, the support of the prior distribution shrinks, as
illustrated in Figure 1.
This explains why values of $K$ close to the recommended value
are preferred.

As an aside, if the true model was unconditional (\textit{i.e.}, without conditional heteroskedasticity),
then there is little loss in using more expanded moments - the extra moments are
superfluous and hence do not change the effective support of the prior. In that case,
no tangible cost in imposed, apart from the computational burden of carrying those moments along.

\begin{table}[h]
\centering
    \begin{tabular}{rcccccccc}
\toprule
      & $Vol(H_{n,K})$ &       & Mean  & SD    & Median & Lower & Upper & Ineff \\
\midrule
\underline{$K=2$} & 0.76  & $\theta_{1}$ & 1.07  & 0.10  & 1.07  & 0.88  & 1.27  & 1.10 \\
      &       & $\theta_{2}$ & 1.01  & 0.14  & 1.01  & 0.74  & 1.29  & 1.14 \\
\midrule
\underline{$K=5$} & 0.73  & $\theta_{1}$ & 1.07  & 0.10  & 1.07  & 0.88  & 1.26  & 1.15 \\
      &       & $\theta_{2}$ & 1.03  & 0.12  & 1.03  & 0.80  & 1.25  & 1.08 \\
\midrule
\underline{$K=10$} & 0.68  & $\theta_{1}$ & 1.07  & 0.09  & 1.07  & 0.89  & 1.25  & 1.14 \\
      &       & $\theta_{2}$ & 1.02  & 0.11  & 1.02  & 0.79  & 1.25  & 1.13 \\
\midrule
\underline{$K=15$} & 0.60  & $\theta_{1}$ & 0.98  & 0.07  & 0.98  & 0.85  & 1.12  & 1.14 \\
      &       & $\theta_{2}$ & 1.10  & 0.10  & 1.10  & 0.91  & 1.29  & 1.13 \\
\midrule
\underline{$K=20$} & 0.54  & $\theta_{1}$ & 0.99  & 0.07  & 0.99  & 0.86  & 1.13  & 1.11 \\
      &       & $\theta_{2}$ & 1.12  & 0.09  & 1.12  & 0.94  & 1.31  & 1.11 \\
\bottomrule
    \end{tabular}
\caption{Example 1 ($n=250$): Volume (prior probability content)
of the convex hull (an estimate of the prior support) and posterior summary for $K$
given by 2, $2n^{1/6}$ and 10, 15, 20.  Results based on 20,000 MCMC draws beyond a
burn-in of 1000. \textquotedblleft Lower\textquotedblright\ and
\textquotedblleft Upper\textquotedblright\ refer to the 0.05 and 0.95
quantiles of the simulated draws, respectively, and \textquotedblleft
Ineff\textquotedblright\ to the inefficiency factor.}\label{Tab:1}
\end{table}

\begin{figure}[h!]

\begin{subfigure}{.33\textwidth}
\begin{center}
\subcaption{$K=5$}
\includegraphics[width=1\textwidth]{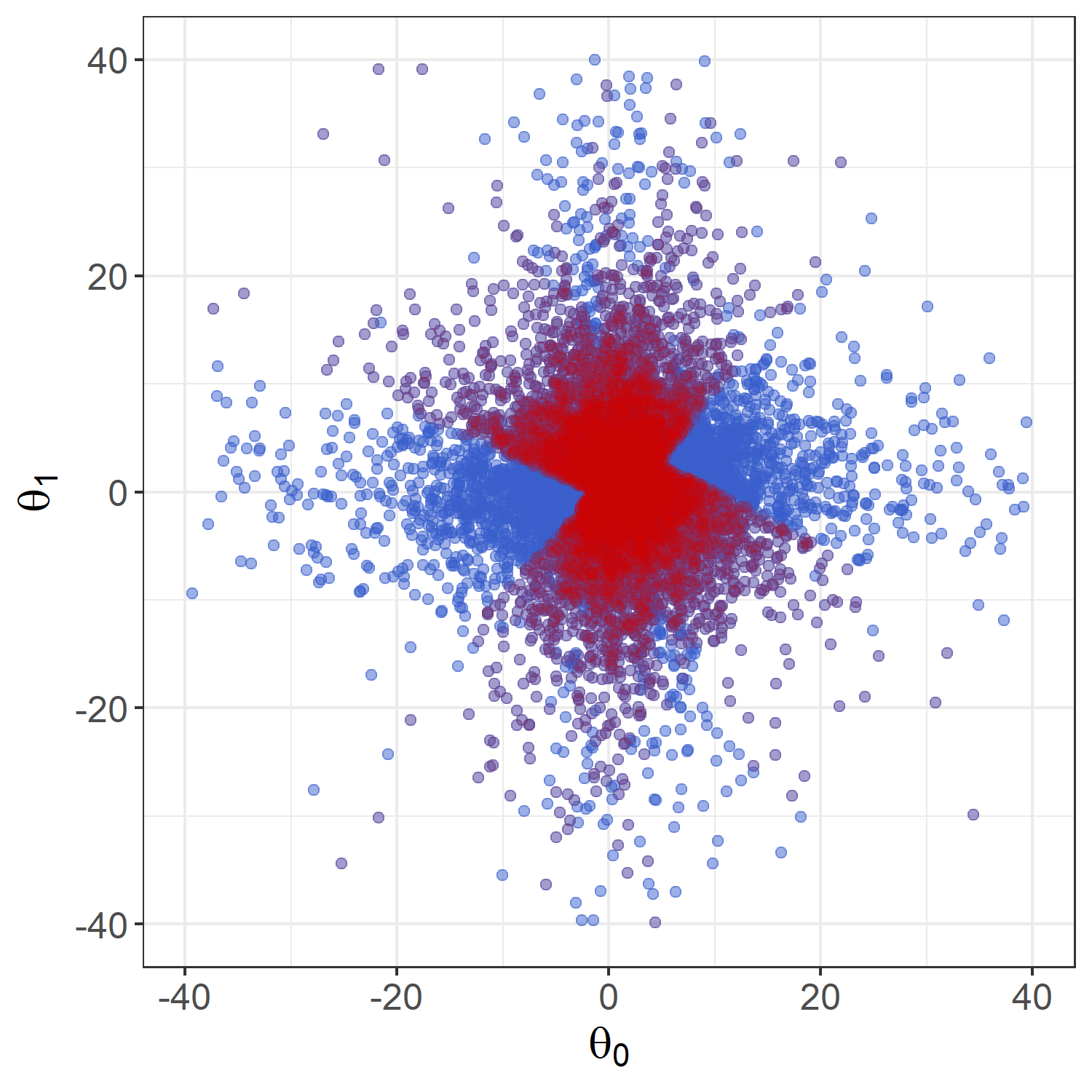}
\end{center}
\end{subfigure}
\begin{subfigure}{.33\textwidth}
\begin{center}
\subcaption{$K=10$}
\includegraphics[width=1\textwidth]{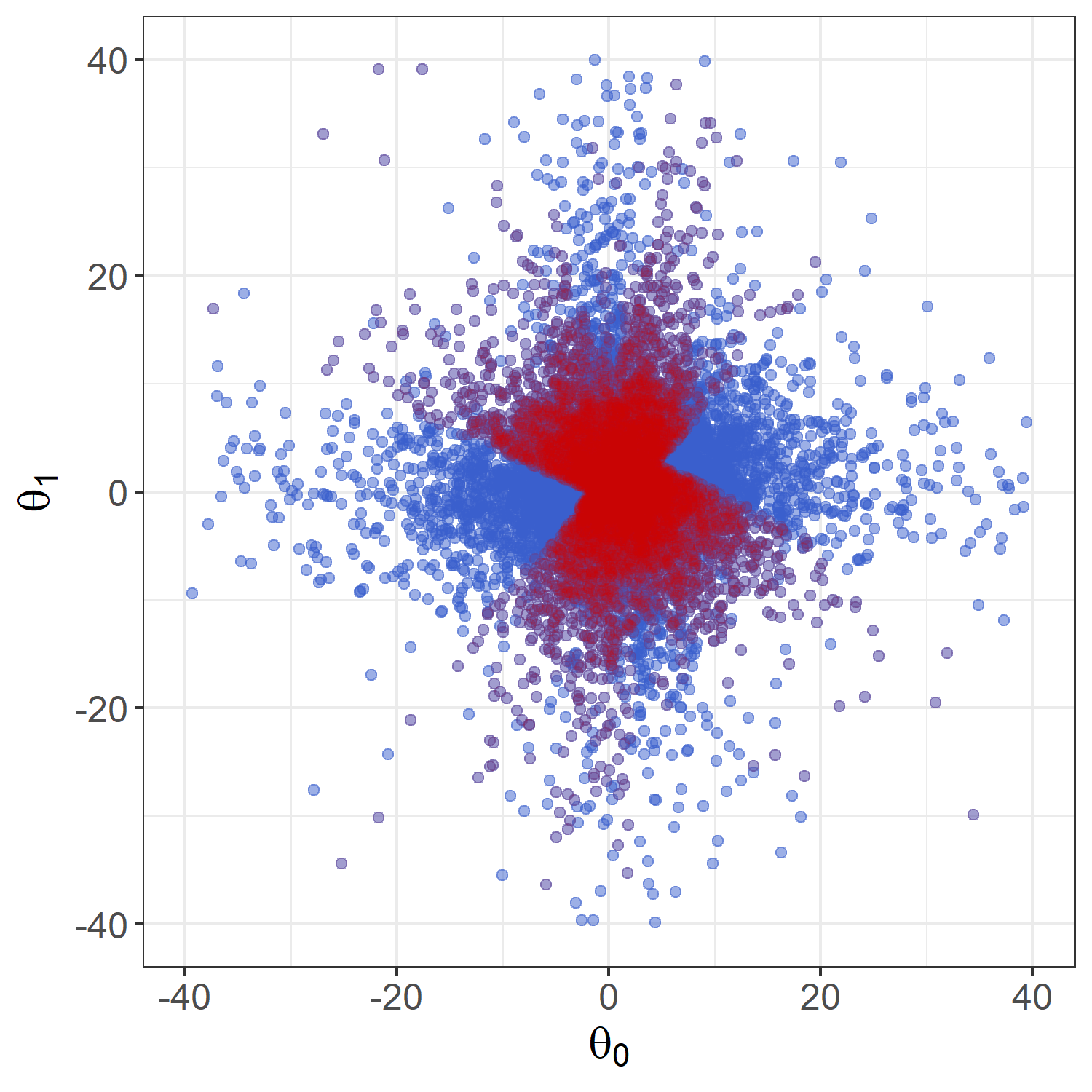}
\end{center}
\end{subfigure}
\begin{subfigure}{.33\textwidth}
\begin{center}
\subcaption{$K=20$}
\includegraphics[width=1\textwidth]{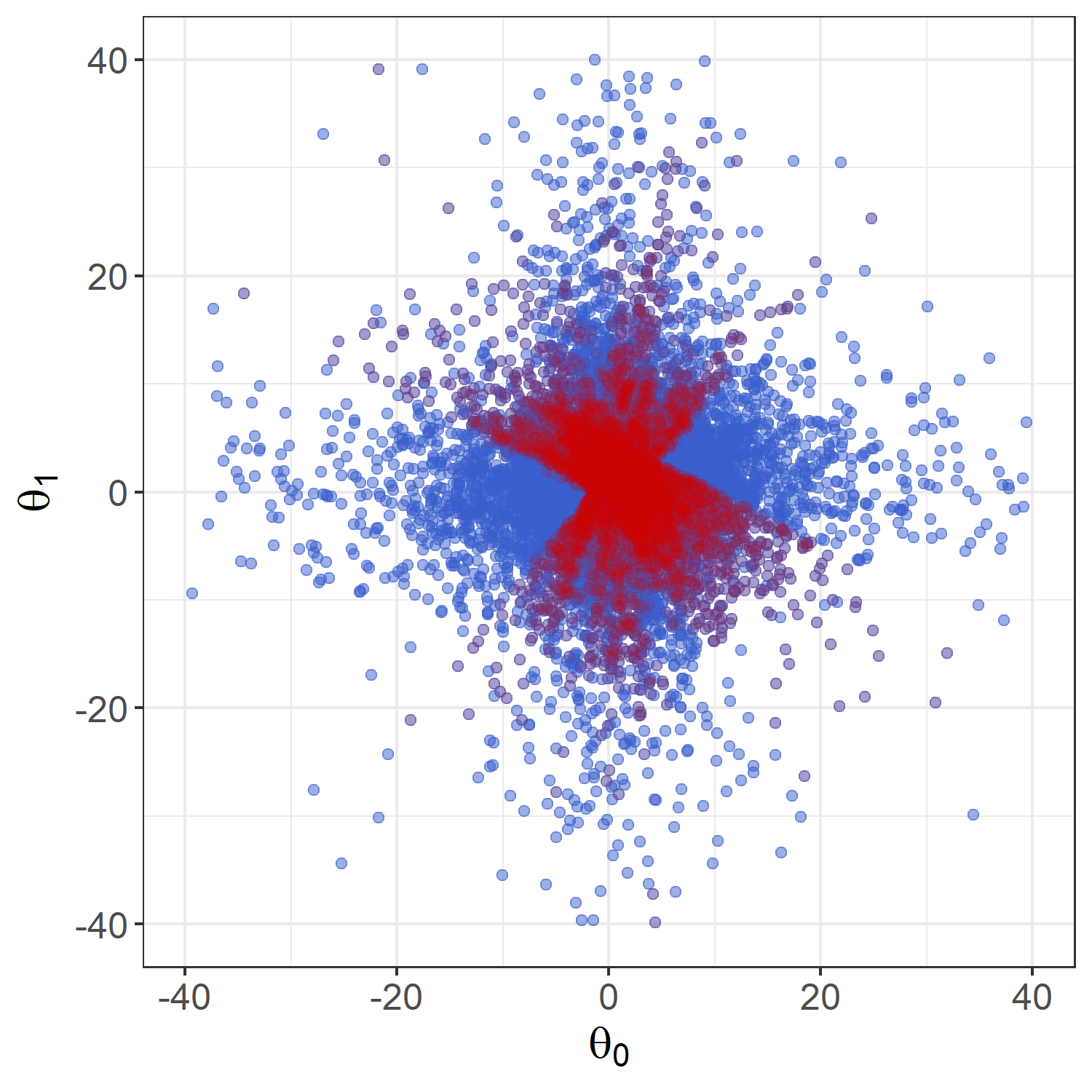}
\end{center}
\end{subfigure}
\caption{Example 1 ($n = 250$): This figure visualizes the volume of the convex hull, $H_{250,K}$,
for different values of $K$. Blue dots are
10,000 independent draws from the default student-t prior. Red dots are the corresponding points in $H_{n,K}$. The volume of $H_{n,K}$ shrinks as $K$ increases (the prior support decreases),
just as shown by the volume estimates presented in Table \ref{Tab:1}. The range in the figures
is set to (-40,40) for visualization clarity.}
\label{fig_ex01_prior}
\end{figure}

\end{contexample1}

\subsection{Asymptotic properties}

\label{ss_Asymptotic_properties} 

Consider now the large sample behavior of the posterior distribution of $\theta$. We let $\theta_*$ and $P_*$, respectively, denote the true value of $\theta$ and of the data distribution $P$. As notation, when the true distribution $P_{*}$ is involved, expectations $\EE^{P}[\cdot ]$ (resp. $\EE^{P}[\cdot |\cdot ]$) are
taken with respect to $P_{*}$ (resp. the conditional distribution associated with $P_{*}$). In addition, we denote $\ell_{n,\theta}(W_i) \coloneqq \log \wh p_{i}(\theta) = \log\frac{e^{\wh{\lambda}(\theta)'g(W_i,\theta)}}{\sum_{j=1}^n e^{\wh{\lambda}(\theta)'g(W_j,\theta)}}$, $$\rho_{\theta}(X,\theta) \coloneqq \frac{\partial \rho(X,\theta)}{\partial\theta^{\prime }},\qquad D(Z) \coloneqq \EE^P[\rho_{\theta}(X,\theta_*)|Z],$$
$$\Sigma(Z)\coloneqq\EE^P[\rho(X,\theta_*)\rho(X,\theta_*)'|Z],\quad \textrm{ and }\quad \rho_{j\theta\theta}(X,\theta_*) \coloneqq \partial^2\rho_j(X,\theta_*)/\partial\theta\partial \theta'.$$
For a vector $a$, $\|a\|$ denotes the Euclidean norm. For a matrix $A$, $\|A\|$ denotes the operator norm (the largest singular value of the matrix). Finally, let $\mathcal{Z}\coloneqq\mathrm{supp}(Z)$ denote the support of $Z$.\\
\indent The first assumption is a normalization for the second moment matrix of the approximating functions which is standard in the literature, see \textit{e.g.} \cite{NEWEY1997} and \cite{DonaldImbensNewey2003}.
\begin{assum}\label{Ass_2_DIN}
For each $K$ there is a constant scalar $\zeta(K)$ such that $\sup_{z\in\mathcal{Z}}\|q^K(z)\| \leq \zeta(K)$, $\EE^P[q^K(Z) q^K(Z)']$ has smallest eigenvalue bounded away from zero
uniformly in $K$, and $\sqrt{K}\leq \zeta(K)$.
\end{assum}
\noindent The bound $\zeta(K)$ is known explicitly in a number of cases depending on the approximating functions we use. \cite{DonaldImbensNewey2003} provide a discussion and explicit formulas for $\zeta(K)$ in the case of splines, power series and Fourier series. We also refer to \cite{NEWEY1997} for primitive conditions for regression splines and power series. 

\begin{assum}\label{Ass_3_DIN}
(a) There exists a unique $\theta_*\in\Theta$ that satisfies $\EE^P[\rho(X,\theta)|Z] = 0$ for the true $P_*$; (b) the data $W_i\coloneqq(X_i,Z_{1i})$, $i=1,\ldots, n$ are i.i.d. according to $P_*$; 
(c) $\EE^P[\sup_{\theta\in\Theta}\|\rho(X,\theta)\|^2|Z]$ is bounded.
\end{assum}

This assumption is the same as \cite[Assumption 3]{DonaldImbensNewey2003}. The following three assumptions are also the same as the ones required by \cite{DonaldImbensNewey2003} to establish asymptotic normality of the Generalized Empirical Likelihood (GEL) estimator.

\begin{assum} \label{Ass_4_DIN}
  (a) $\theta_*\in int(\Theta)$; (b) $\rho(X,\theta)$ is twice continuously differentiable in a neighborhood $\mathcal{U}$ of $\theta_*$, $\EE^P[\sup_{\theta\in\mathcal{U}}\|\rho_{\theta}(X,\theta)\|^2|Z]$ and $\EE^P[\sup_{\theta\in\mathcal{U}}\|\rho_{j\theta\theta}(X,\theta_*)\|^2|Z]$, $j=1,\ldots d$, are bounded on $\mathcal{Z}$; (c) $\EE^P[D(X)D(X)']$ is nonsingular.
\end{assum}

\begin{assum}\label{Ass_5_DIN}
(a) $\Sigma(Z)$ has smallest eigenvalue bounded away from zero; (b) for a neighborhood $\mathcal{U}$ of $\theta_*$, $\EE^P[\sup_{\theta\in\mathcal{U}}\|\rho(X,\theta)\|^4|z]$ is bounded, and for all $\theta\in\mathcal{U}$, $\|\rho(X,\theta) - \rho(X,\theta_*)\| \leq
\delta(X)\|\theta - \theta_*\|$ and $\EE^P[\delta(X)^2|Z]$ is bounded.
\end{assum}

\begin{assum}\label{Ass_6_DIN}
  There is $\gamma > 2$ such that $\EE^P[\sup_{\theta\in\Theta}\|\rho(X,\theta)\|^{\gamma}] < \infty$ and $\zeta(K)^2 K/ n^{1 - 2/\gamma} \rightarrow 0$.
\end{assum}

Part (b) of Assumption \ref{Ass_5_DIN} imposes a Lipschitz condition which allows application of uniform convergence results. The last assumption is about the prior distribution of $\theta$ and is standard in the Bayesian literature on frequentist asymptotic properties of Bayes procedures.
\begin{assum}\label{Ass_prior}
  (a) $\pi$ is a continuous probability measure that admits a density with respect to the Lebesgue measure; (b) $\pi$ is positive on a neighborhood of $\theta_*$.
\end{assum}

We are now able to state our first major result in which we establish the asymptotic normality and efficiency of the posterior distribution of the local parameter $h \coloneqq \sqrt{n}(\theta - \theta_*)$.

\begin{thm}[Bernstein-von Mises]\label{thm_BvM_correctly_specified}
 Under Assumptions \ref{Ass_1_DIN}-\ref{Ass_prior}, if $K\rightarrow \infty$, $\zeta(K) K^2/\sqrt{n} \rightarrow 0$, 
 and if for any $\delta > 0$, $\exists\epsilon>0$ such that as $n\rightarrow
\infty$
\begin{equation}\label{Ass_identification_rate_contraction_correct_specification}
P\left(\sup_{\|\theta - \theta_*\| > \delta}\frac{1}{n}\sum_{i=1}^n \left(\ell_{n,\theta}(W_i) - \ell_{n,\theta_*}(W_i)\right) \leq - \epsilon\right) \rightarrow 1,
\end{equation}
then the posterior distribution $\pi(\sqrt{n}(\theta - \theta_*)|W_{1:n})$ converges in total variation towards a random Normal distribution, that is,
\begin{equation}\label{BVM _result_correct_specification}
  \sup_{B}\left|\pi(\sqrt{n}(\theta-\theta_*)\in B|W_{1:n},K) - \mathcal{N}_{\Delta_{n,\theta_*},V_{\theta_*}}(B)\right|\overset{p}{\to} 0,
\end{equation}
where $B\subseteq \Theta$ is any Borel set, $\Delta_{n,\theta_*} \coloneqq -\frac{1}{\sqrt{n}}\sum_{i=1}^n V_{\theta_*}D(Z_i)'\Sigma(Z_i)^{-1}\rho(X_i,\theta_*)$ is bounded in probability and $V_{\theta_*} \coloneqq \left(\EE^P[D(Z)' \Sigma(Z)^{-1} D(Z)]\right)^{-1}$.
\end{thm}
We note that the centering $\Delta_{n,\theta_*}$ of the limiting normal distribution satisfies $\frac{1}{\sqrt{n}}\sum_{i=1}^n\frac{d\log\wh p_i(\theta_*)}{d\theta} - V_{\theta_*}^{-1}\Delta_{n,\theta_*} \overset{p}{\to} 0$. We also note that the condition $\zeta(K) K^2/\sqrt{n} \rightarrow 0$ in the theorem implies $K/n\rightarrow 0$, which is a classical condition in the sieve literature. This condition is required to establish a stochastic Local Asymptotic Normality (LAN) expansion, which is an intermediate step to prove the BvM result, as we explain below. The LAN expansion is not required to establish asymptotic normality of the GEL estimators, which explains why our condition is slightly stronger than the condition $\zeta(K) K/\sqrt{n} \rightarrow 0$ required by \citet*{DonaldImbensNewey2003}. 
On the other hand, our condition is weaker than the condition $\zeta(K)^2 K^2/\sqrt{n} \rightarrow 0$ required by \citet*{DONALD2009} to establish the mean square error of the GEL estimators. The asymptotic covariance of the posterior distribution coincides with the semiparametric efficiency bound given in \cite{Chamberlain1987} for conditional moment condition models. This means that, for every $\alpha\in (0,1)$, $(1 - \alpha)$-credible regions constructed from the posterior of $\theta$ are $(1 - \alpha)$-confidence sets asymptotically.

The proof of this theorem is given in the supplementary appendix and consists of three steps. In the first step we show consistency of the posterior distribution of $\theta$, namely:
\begin{equation}\label{eq_posterior_consistency_manuscript}
  \pi\left(\left.\sqrt{n}\|\theta - \theta_*\| > M_n\right|W_{1:n},K\right)\overset{p}{\to} 0
\end{equation}
\noindent for any $M_n\rightarrow \infty$, as $n\rightarrow \infty$. To show this, the identification assumption \eqref{Ass_identification_rate_contraction_correct_specification} is used. In the second step we show that the ETEL function satisfies a stochastic LAN expansion:
  \begin{equation}\label{eq_stochastic_LAN_manuscript}
    \sup_{h\in \mathcal{H}}\left|\sum_{i=1}^n \ell_{n,\theta_* + h/\sqrt{n}}(W_i) - \sum_{i=1}^n \ell_{n,\theta_*}(W_i) - h' V_{\theta_*}^{-1}\Delta_{n,\theta_*} + \frac{1}{2}h'V_{\theta_*}^{-1}h\right| = o_p(1),
  \end{equation}
where $\mathcal{H}$ denotes a compact subset of $\mathbb{R}^p$ and $V_{\theta_*}^{-1}\Delta_{n,\theta_*}\overset{d}{\to} \mathcal{N}(0, V_{\theta_*}^{-1})$. As the ETEL function is an integrated likelihood, expansion \eqref{eq_stochastic_LAN_manuscript} is better known as integral LAN in the semiparametric Bayesian literature, see \textit{e.g.} \cite[Section 4]{BickelKleijn2012}. 
In the third step of the proof we use arguments as in the proof of \cite[Theorem 10.1]{VanDerVaart2000} to show that \eqref{eq_posterior_consistency_manuscript} and \eqref{eq_stochastic_LAN_manuscript} imply asymptotic normality of $\pi(\sqrt{n}(\theta-\theta_*)\in B|W_{1:n},K)$. While these three steps are classical in proving the Bernstein-von Mises phenomenon, establishing \eqref{eq_stochastic_LAN_manuscript} raises challenges that are otherwise absent. This is because the ETEL function is a nonstandard likelihood that involves estimated parameters $\wh\lambda(\theta_*)$ whose dimension is $d K$, which increases with $n$. 
While $\|\wh\lambda(\theta_*)\|$ and $\|\frac{1}{n}\sum_{i=1}^n g(W_i,\theta_*)\|$ are expected to converge to zero in the correctly specified case, the rate of convergence is slower than $n^{-1/2}$. In the supplementary appendix we show that this rate is $\sqrt{K/n}$ under the previous assumptions.

\subsection{Misspecified model}\label{s_Asymptotic_result_ms}

We now generalize the preceding BvM result for the important class of
misspecified conditional moment models.

\begin{df}[Misspecified model]\label{def_misspecified_model}
  We say that the conditional moment conditions model $\E^P[\rho(X,\theta)|Z] = 0$ is misspecified if the set of probability measures implied by the moment restrictions does not contain the true data generating process $P_*$ for any $\theta\in \Theta$, that is, $P_*\notin \mathcal{P}$ where $\mathcal{P} \coloneqq \bigcup_{\theta\in\Theta }\wtl{\mathcal{P}}_{\theta}$ and $\wtl{\mathcal{P}}_{\theta} = \{Q\in \mathbb{M}_{X|Z}; \,\EE^Q[\rho(X,\theta)|Z] = 0 \; a.s.\}$ with $\mathbb{M}_{X|Z}$ the set of all conditional probability measures of $X|Z$.
\end{df}
\indent In essence, if \eqref{eq:condmom} is misspecified then there is no $\theta\in\Theta$ such that $\EE^P[\rho(X,\theta)\otimes q^K(Z)] = 0$ almost surely for every $K$ large enough. Now, for every $\theta\in\Theta$ define $Q^*(\theta)$ as the minimizer of the Kullback-Leibler divergence of $P_*$ to the model $\mathcal{P}_{\theta} \coloneqq \{Q\in \mathbb{M}; \EE^Q[g(W,\theta)]=0\}$, where $\mathbb{M}$ denotes the set of all the probability measures on $\mathbb{R}^{d_w}$. That is, $Q^*(\theta) \coloneqq \mathrm{arginf}_{Q\in\mathcal{P}_{\theta}}\mathbb{K}(Q||P_*)$, where $\mathbb{K}(Q||P_*)\coloneqq\int\log(dQ/dP_*)dQ$. If we suppose that the dual representation of the Kullback-Leibler minimization problem holds, then the $P_*$-density of $Q^*(\theta)$ has the closed form: $[dQ^*(\theta)/dP_*](W_i) = \frac{e^{\lambda_{\circ}'g(W_i,\theta)}}{\EE^P[ e^{\lambda_{\circ}'g(W_j,\theta)}]}$, where $\lambda_\circ$ denotes the tilting parameter and is defined in the same way as in the correctly specified case:
\begin{equation}\label{eq_def_tilting_parameter_miss}
  \lambda_{\circ} \coloneqq \lambda_{\circ}(\theta)\coloneqq\arg\min_{\lambda\in\mathbb{R}^{d K}}\EE^P[e^{\lambda'g(W_i,\theta)}].
\end{equation}

We also impose a condition to ensure that the probability measures $\mathcal{P}\coloneqq\bigcup_{\theta\in\theta}\mathcal{P}_\theta$, which are implied by the model, are dominated by the true probability measure $P_*$. This is required for the validity of the dual theorem.
Therefore, following \cite[Theorem 3.1]{Sueishi2013}, we replace Assumption \ref{Ass_3_DIN} (a) by the following.
\begin{assum}\label{Ass_absolute continuity}
  For every $\theta\in\Theta$, there exists $Q\in \mathcal{P}_{\theta}$ such that $Q$ is mutually absolutely continuous with respect to $P_*$, where $\mathcal{P}_{\theta}\coloneqq\{Q\in \mathbb{M}; \EE^Q[g(W,\theta)]=0\}$ and $\mathbb{M}$ denotes the set of all the probability measures on $\mathbb{R}^{d_w}$.
\end{assum}
This assumption implies that $\mathcal{P}_{\theta}$ is non-empty. A similar assumption is also made by \citet{kleijn2012} and \citet*{ChibShinSimoni2018} to establish the BvM under misspecification. The pseudo-true value of the parameter $\theta\in \Theta$ is denoted by $\theta_\circ$ and is defined as the minimizer of the Kullback-Leibler divergence between the true $P_*$ and $Q^*(\theta)$:
\begin{equation}\label{eq_pseudo_true_value}
\theta_\circ \coloneqq \mathrm{arginf}_{\theta\in \Theta}\mathbb{K}(P_*||Q^*(\theta)),
\end{equation}
where $\mathbb{K}(P_*||Q^*(\theta))\coloneqq\int \log(dP_*/dQ^*(\theta)) dP_*$. Under the preceding absolute continuity assumption, the pseudo-true value $\theta_\circ$ is available as
\begin{equation}\label{eq_def_pseudo_true_value}
  \theta_\circ = \argmax_{\theta\in\Theta}\EE^P\log\left(\frac{e^{\lambda_{\circ}'g(W_i,\theta)}}{\EE^P[ e^{\lambda_{\circ}'g(W_j,\theta)}]}\right).
\end{equation}
\noindent Note that $\lambda_\circ (\theta_\circ)$, the value of the tilting parameter at the pseudo-true value $\theta_\circ$, is nonzero because the moment conditions do not hold.\\
\indent Assumption \ref{Ass_absolute continuity} implies that $\mathbb{K}(Q^*(\theta_\circ)||P_*)<\infty$. We supplement this with the assumption that $\mathbb{K}(P_*||Q^*(\theta))<\infty$, $\forall \theta\in\Theta$ (so that $\mathbb{K}(P_*||Q^*(\theta_\circ))<\infty$). Because consistency in misspecified models is defined with respect to the pseudo-true value $\theta_\circ$, we need to replace Assumption \ref{Ass_prior} \textit{(b)} by the following Assumption \ref{Ass_prior_BvM_misspecified} (b) which, together with Assumption \ref{Ass_prior_BvM_misspecified} \textit{(a)}, requires the prior to put enough mass to balls around $\theta_\circ$.
\begin{assum}\label{Ass_prior_BvM_misspecified}
  (a) $\pi$ is a continuous probability measure that admits a density with respect to the Lebesgue measure; (b) The prior distribution $\pi$ is positive on a neighborhood of $\theta_\circ$, where $\theta_\circ$ is as defined in \eqref{eq_def_pseudo_true_value}.
\end{assum}

Hereafter, we use the sub/super index $Q^*(\theta_\circ)$ to denote an expectation, a variance or covariance taken with respect to the probability $Q^*(\theta_\circ)$. The following assumption is analogous to the second part of Assumption \ref{Ass_2_DIN} for the $P_*$-density of $Q_*(\theta)$ replacing $P$.
\begin{assum}\label{Ass_2_DIN_Miss}
For each $K$ the matrix $\EE^{Q^*(\theta_\circ)}[q^K(Z) q^K(Z)']$ has smallest eigenvalue bounded away from zero uniformly in $K$.
\end{assum}

In the next assumption we denote by $int(\Theta)$ the interior of $\Theta$ and by $\mathcal{U}$ a ball centered at $\theta_\circ$ with radius $h/\sqrt{n}$ for some $h\in\mathcal{H}$ and $\mathcal{H}$ a compact subset of $\mathbb{R}^p$.
\begin{assum}\label{Ass_3_BvM_misspecified}
  (a) The data $W_i\coloneqq(X_i,Z_i)$, $i=1,\ldots, n$ are i.i.d. according to $P_*$ and \\
  (b) The pseudo-true value $\theta_\circ\in int(\Theta)$ is the unique maximizer of $$\lambda_\circ(\theta)'\EE^P[g(W,\theta)] - \log\EE^P[\exp\{\lambda_\circ(\theta)'g(W,\theta)\}],$$ where $\Theta\subset\mathbb{R}^p$;\\
  (c) $\rho(X,\theta)$ is continuous at each $\theta\in\Theta$ with probability one;\\
  (d) $\rho(X,\theta)$ is twice continuously differentiable in a neighborhood $\mathcal{U}$ of $\theta_\circ$ with probability one and for $\kappa = 0,1$, 
  $\EE^P\left[\sup_{\theta\in\mathcal{U}}\left|\frac{dQ^*(\theta)}{dP_*}(W)\right|^{\kappa}\|\rho_{j\theta\theta}(X,\theta)\|^2\|q^K(Z)\|^2\right] = \mathcal{O}(K)$, $j=1,\ldots d$;\\
  (e) for a neighborhood $\mathcal{U}$ of $\theta_\circ$ and for $\kappa =0,1,2$, $j=2,4$ it holds that
  $$\EE^P\left[\sup_{\theta\in\mathcal{U}}\left|\frac{dQ^*(\theta)}{dP_*}(W_i)\right|^{\kappa}\|\rho(X,\theta)\|^j\|q^K(Z)\|^j\right] = \mathcal{O}(\zeta(K)^{j-2} K),$$ where $\zeta(K)$ is as defined in Assumption \ref{Ass_2_DIN};\\
  (f) for a neighborhood $\mathcal{U}$ of $\theta_\circ$ and for $\kappa = 0,1,2$, $j=1,2,4$ it holds that
  $$\EE^P\left[\sup_{\theta\in\mathcal{U}}\left|\frac{dQ^*(\theta)}{dP_*}(W_i)\right|^{\kappa}\|\rho_{\theta}(X_i,\theta)\|^j\|q^K(Z)\|^j\right] = \mathcal{O}(\zeta(K)^{\max\{j-2,0\}} K),$$ 
  where $\zeta(K)$ is as defined in Assumption \ref{Ass_2_DIN};\\
  (g) the matrix $\EE^{Q^*(\theta_\circ)}\left[\left.\rho(X,\theta_\circ)\rho(X,\theta_\circ)'\right|Z\right]$ has smallest (resp. largest) eigenvalue bounded away from zero (resp. infinity);\\
  (h) for a neighborhood $\mathcal{U}$ of $\theta_\circ$ it holds that $\EE^P\left[\sup_{\theta\in\mathcal{U}}\left|\frac{dQ^*(\theta)}{dP_*}(W_i)\right|^{2}\right]$ is bounded. 
\end{assum}
Assumption \ref{Ass_3_BvM_misspecified} (b) guarantees uniqueness of the pseudo-true value and is a standard assumption in the literature on misspecified models (see \textit{e.g.} \cite{White1984}). 
Assumptions \ref{Ass_3_BvM_misspecified} (d)-(f) are the counterparts of Assumptions \ref{Ass_4_DIN} (b) and \ref{Ass_5_DIN} (b), respectively, for the misspecified case. It is important to notice that they implicitly contain the first part of Assumption \ref{Ass_2_DIN}. The reason why we cannot separate the part involving the moment function $\rho(X,\theta)$ (or its derivative) and the one involving $q^K(Z)$ in the assumption, as we do for the correctly specified model, is that the $P_*$-density of $Q_*(\theta)$ cannot be factorized in a conditional density of $X$ given $(Z,\theta)$ and a marginal density of $Z$ independent of $\theta$. In particular, in the misspecified case the pseudo-true value of the tilting parameter $\lambda_\circ(\theta_\circ)$ is not equal to zero as it is the tilting parameter in the correctly specified case. Assumption \ref{Ass_3_BvM_misspecified} (g) is the counterpart of Assumption \ref{Ass_5_DIN} (b) for the misspecified case.
%
%

The BvM theorem for misspecified models now follows. Let $G_i(\theta_\circ) \coloneqq \rho_{\theta}(X_i,\theta_\circ)\otimes q^K(Z_i)$, $D^{\dagger}_\circ(Z) \coloneqq \EE^{Q^*(\theta_\circ)}[\rho_{\theta}(X,\theta_\circ)|Z]$ and $\Sigma_{\circ}(Z)\coloneqq\EE^{Q^*(\theta_\circ)}[\rho(X,\theta_\circ)\rho(X,\theta_\circ)'|Z]$. Moreover, let $\mathcal{H}$ denote a compact subset of $\mathbb{R}^p$ and $\theta_h \coloneqq \theta_\circ + h/\sqrt{n}$, with $h\in\mathcal{H}$.
\begin{thm}[Bernstein-von Mises (misspecified)]\label{thm_BvM_misspecified}
 Let Assumptions \ref{Ass_1_DIN}, \ref{Ass_2_DIN}, \ref{Ass_absolute continuity} - \ref{Ass_3_BvM_misspecified} hold. Assume that there exists a constant $C>0$ such that for any sequence $M_n\rightarrow \infty$,
  \begin{equation}\label{Ass_identification_rate_contraction_appendix}
    P_*\left(\sup_{\|\theta - \theta_\circ\| > M_n/\sqrt{n}}\frac{1}{n}\sum_{i=1}^n \left(\ell_{n,\theta}(W_i) - \ell_{n,\theta_\circ}(W_i)\right) \leq - CM_n^2/n\right) \rightarrow 1,
  \end{equation}
  \noindent as $n\rightarrow \infty$. 
  If $K\rightarrow \infty$ and $\sup_{\theta\in\mathcal{U}}\|\lambda_\circ(\theta)\|^2\max\{\zeta(K),K\}K\sqrt{K/n} \rightarrow 0$ then, the posteriors converge in total variation towards a Normal distribution, that is,
\begin{equation}\label{BVM _result_misspecified_specification}
  \sup_{B}\left|\pi(\sqrt{n}(\theta-\theta_\circ)\in B|W_{1:n},K) - \mathcal{N}_{\Delta_{n,\theta_\circ},V_{\theta_\circ}}(B)\right|\overset{p}{\to} 0,
\end{equation}
where $B\subseteq \Theta$ is any Borel set, $\Delta_{n,\theta_\circ}$ is a random vector bounded in probability and $V_{\theta_\circ}$ is a positive definite matrix 
 equal to the inverse of:
{\footnotesize
\begin{multline*}
  V_{\theta_\circ}^{-1} = \EE^P\left[D^{\dagger}_\circ(Z)\Sigma_{\circ}(Z)^{-1}D_\circ(Z)\right] + \EE^{Q^*(\theta_\circ)} \left[G_i(\theta_\circ)'\lambda_\circ(\theta_\circ) g(W_i,\theta_\circ)'\right]\Omega_\circ^{-1}\EE^P[G_i(\theta_\circ)]\\
  - \frac{d\lambda_{\circ}(\theta_\circ)'}{d\theta}\left(\EE^{P}\left[G_i(\theta_{\circ})\right] - \EE^{Q^*(\theta_\circ)}\left[G_i(\theta_{\circ})\right]\right) - \sum_{j=1}^{d} \frac{d^2\lambda_{\circ,j}(\theta_\circ)'}{d\theta d\theta'}\EE^P[\rho_{j}(X_i,\theta_\circ)q^K(Z)]\\
  - \sum_{j=1}^{d} \left(\EE^P\left[\rho_{j\theta\theta}(X_i,\theta_\circ)q^K(Z_i)\right] - \EE^{Q^*(\theta_\circ)}\left[\rho_{j\theta\theta}(X_i,\theta_\circ)q^K(Z_i)\right]\right) \lambda_{\circ,j}(\theta_\circ)\\
  + Var_{Q^*(\theta_\circ)}[G_i(\theta_\circ)'\lambda_\circ(\theta_\circ)] + \frac{d\lambda_\circ(\theta_\circ)'}{d\theta}Cov_{Q^*(\theta_\circ)}\left(g(W_i,\theta_\circ),G_i(\theta_\circ)'\lambda_\circ(\theta_\circ)\right).
\end{multline*}
}
\end{thm}
Just as in \cite{kleijn2012}, this theorem establishes that the posterior distribution of the centered and scaled parameter $\sqrt{n}(\theta-\theta_\circ)$ converges to a Normal distribution with a random mean that is bounded in probability. The rate restriction $\sup_{\theta\in\mathcal{U}}\|\lambda_\circ(\theta)\|^2\max\{\zeta(K),K\}K\sqrt{K/n} \rightarrow 0$ is much stronger than the one in Theorem \ref{thm_BvM_correctly_specified}. The slower rate condition on $K$ is intuitive. When the conditional moment conditions are misspecified, limiting the number of implied unconditional moments serves to limit the magnification of the misspecification. In the (completely) hypothetical situation where one knew that the conditional moment conditions are misspecified, one would either discard the misspecified moment conditions or take a small and fixed number of expanded moment conditions (for instance with $K=1$). In practice, of course, this strategy cannot be implemented (because one does not know whether the model is correctly specified or misspecified) and, therefore, $K$ must always go to $\infty$, but slower than under correct specification.

An additional remark is that, because of misspecification, the covariance matrix $V_{\theta_\circ}$ appearing in Theorem \ref{thm_BvM_misspecified} is expected to be different from the asymptotic covariance matrix of the corresponding frequentist point estimator in \cite[Theorem 4.1]{AI2007}. Therefore, while correctly centered, the $(1 - \alpha)$ posterior credible sets are not in general $(1 - \alpha)$ confidence sets. The coverage rate can be less, or more, than the nominal level depending on the true data generating process and the extent of misspecification.

The strategy of the proof of this result is generally similar to the proof of Theorem \ref{thm_BvM_correctly_specified} with $\theta_*$ replaced by the pseudo-true value $\theta_\circ$. However, proving that the ETEL function satisfies a stochastic LAN expansion is more complex, for the following reasons.
First, the limit of $\wh\lambda(\theta_\circ)$ is $\lambda_\circ(\theta_\circ)$, which is not zero. Therefore, several terms that were equal to zero in the LAN expansion under correct specification, are non-zero in the misspecified case. The limit in distribution of these terms has to be derived. This explains our stronger assumptions with respect to the correctly specified case. Second, the quantity $\frac{1}{\sqrt{n}}\sum_{i=1}^n g(W_i,\theta_\circ)$ is no longer centered on zero, which leads to an additional bias term.
%


\begin{contexample1}

 Consider now the impact of $K$ in misspecified models. For simplicity, we fix $\theta_{0}$ to a certain value and we only estimate $\theta_{1}$. We generate the sample data as above with $(\theta_{0}=1,\theta_1 = 1)$ and then estimate $\theta_1$, fixing $\theta_{0}$ at 0.5. Therefore,
 the moment condition $\EE^{P}[\varepsilon |Z]=0$ is misspecified. The pseudo-true value of $\theta_{1}$ is obtained from a side calculation. More specifically, we generate 5 million observations from the true model. We assume that this sample size is large enough to represent the $n \rightarrow \infty$ situation. Then, we misspecify the CM conditions by setting $\theta_{0}=0.5$, construct the expanded moments with these 5 millions observations setting $K = 26$, and numerically find the value of $\theta_{1}$ that maximizes the BETEL posterior distribution. The pseudo true value is this maximized value. It is 1.004. The adverse impact of increasing $K$ (for a fixed $n$) on the Bayesian bias and posterior standard deviation is reported in Table \ref{tab:ex01mis freq}. As in the case of the correctly specified models, relatively small values of $K$ (around $K=5$ for $n=250$) lead to the best results and, in addition, the value of the Bayesian bias increase when $K$ increases beyond the recommended value. The posterior sd declines more sharply as K increases. As pointed out before, these effects are due to the reduction in the support of the prior, now magnified by moment misspecification. Finally, we calculate frequentist coverage rate of the equal-tailed 90\% credible set based on 100 repetitions for both the correctly specified model ($\theta_0$ now set equal to 1) and misspecified model with $K\approx2n^{1/6}$. Consistent with our theory the BETEL credible set is different from the  frequentist confidence set when conditional moment conditions are misspecified. When $n=250$ the coverage rates are 91\% for the correctly specified case and 86\% for the misspecified case. When the number of observations increases to 1,000, the coverage rate for the misspecified case further moves down to 80\% while the coverage rate for the correctly specified case remains around 90\%.

\begin{table}[h]
\centering
    \begin{tabular}{l|cc|cccc}
    \toprule
          & \multicolumn{2}{c|}{Correctly specified model} & \multicolumn{2}{c}{Misspecified model} \\
    $n = 250$ & $|\mathrm{Bias}|$ & $\mathrm{SD}$ & $|\mathrm{Bias}|$ & $\mathrm{SD}$   \\
    \midrule
$K=2$ & 0.277 & 0.147 & 0.273 & 0.154 \\
$K=5$ & 0.048 & 0.110 & 0.044 & 0.102 \\
$K=9$ & 0.053 & 0.105 & 0.049 & 0.098 \\
$K=12$ & 0.059 & 0.105 & 0.055 & 0.097 \\
$K=20$ & 0.064 & 0.102 & 0.060 & 0.084 \\
    \midrule
    \bottomrule
    \end{tabular}%
\centering
\caption{Bayesian bias and posterior sd for different values of $K$ under correct and incorrect conditional moments.}
\label{tab:ex01mis freq}
\end{table}

\end{contexample1}

\section{Model Comparisons}

\label{s_4}

\label{s_Model_Comparison}

In practice, we can be unsure about elements of
the conditional moment model. For instance, we can be faced
with a large number of variables in $Z$, but only some of
which are relevant. In such cases, any specific model may be
considered to be misspecified, and the goal is
to find the best model given the data.

Let $M_{\ell }$ denote the $\ell $th model in the model space.
Each model is characterized by a parameter $\theta ^{\ell }$ and an extended
set of moment functions given by $g^{\ell }(W,\theta ^{\ell })$.
In addition, each model $M_{\ell }$ is described by a
prior distribution for $\theta ^{\ell }$. The posterior distribution
is obtained based on \eqref{eq_betel posterior}. The aim is to compare these
models by marginal likelihoods, denoted by $m(W_{1:n}|M_{\ell},K)$. These are each calculated by the marginal likelihood identity of \cite{Chib1995} (where we explicit the dependence on $M_{\ell}$ in the notation):
\begin{equation}
  \log m(W_{1:n}|M_{\ell },K)=\log \pi (\tilde{\theta}^{\ell }|M_{\ell })+\log p(W_{1:n}|\tilde{\theta}^{\ell },M_{\ell },K)-\log \pi (\tilde{\theta}^{\ell}|W_{1:n},M_{\ell }, K),  \label{eq_MSC_ML}
\end{equation}%
and by the method of \cite{ChibJeliazkov2001}. In this expression, $\tilde{\theta}^{\ell }$ is any point in the support of the posterior (such as the posterior mean).

\begin{rem}
  Comparison of CM condition models differs in one important
aspect from the framework for comparing unconditional moment condition
models that was established in \citet*{ChibShinSimoni2018}, where it is shown
that to make the unconditional moment condition models comparable it is
necessary to linearly transform the moment functions so that all the
transformed moments are included in each model. This linear transformation
consists of adding an extra parameter different from zero to the components
of the vector $g(\theta ,W)$ that correspond to the restrictions not
included in a specific model. When comparing conditional moment models,
however, this transformation is not necessary because the convex hulls
associated with different expanded models have the same dimension
asymptotically.
\end{rem}

\subsection{Model selection consistency}

Suppose that there are $J$ contending models. Suppose also that at least $J-1$ of these models are misspecified and the remaining one can be either misspecified or correctly specified. Moreover, suppose that a model $M_\ell$ is selected by the size of the marginal likelihoods. Then, in Theorem \ref{thm_consistency_selection_misspecification} we show that this criterion in the limit picks the model $M_\ell$ with the smallest Kullback-Leibler divergence between $P_*$ and the corresponding $Q^*(\theta^\ell)$, where $Q^*(\theta^\ell) = \arginf_{Q\in\mathcal{P}_{\theta^{\ell}}}\mathbb{K}(Q||P_*)$ and $\mathcal{P}_{\theta^{\ell}}$ is defined in
Section \ref{s_Asymptotic_result_ms}.

\begin{thm}
\label{thm_consistency_selection_misspecification} Let the assumptions of
Theorem \ref{thm_BvM_misspecified} hold. Let us consider the comparison of $J<\infty$ models $M_\ell$, $\ell=1,\ldots,J$, such that $J-1$ of these models each has at least one misspecified moment condition and model $M_j$ can be either correctly specified or contain some misspecified moment condition, that is, $M_{\ell}$ does not satisfy Assumption \ref{Ass_3_DIN} (a), $\forall \ell \neq j$. Then,
\begin{equation*}
\lim_{n\rightarrow \infty}P_*\left(\log m(W_{1:n}|M_{j},K) > \max_{\ell\neq j}\log m(W_{1:n}|M_{\ell}.K)\right) = 1
\end{equation*}
if and only if $\mathbb{K}(P_*||Q^*(\theta_{\circ}^j))< \min_{\ell\neq j}\mathbb{K}(P_*||Q^*(\theta_{\circ}^{\ell}))$, where $\mathbb{K}(P||Q):=\int \log(dP/dQ)dP$.
\end{thm}

Note that if one model in the contending set of models is correctly specified, then this model will have zero Kullback-Leibler divergence and, therefore, according to Theorem \ref{thm_consistency_selection_misspecification}, that model will have the largest marginal likelihood and will be selected by our procedure.

To understand the ramifications of the preceding result, suppose that we are
interested in comparing models with the same moment conditions but different
conditioning variables:
\begin{equation}
\text{Model 1: }\mathbf{E}^{P}[\rho (X,\theta )|Z_{1}]=0,\qquad \text{Model
2: }\mathbf{E}^{P}[\rho (X,\theta )|Z_{2}]=0,
\label{eq_general_case_3rd_situation}
\end{equation}%
where $Z_{1}$ and $Z_{2}$ may have some elements in common, in particular $%
Z_{2}$ might be a subvector of $Z_{1}$ (or vice versa). A situation of this
type, where we are unsure about the validity of instrumental variables, is
the following.

\begin{example2} \label{example2}
(Comparing IV models) Consider the following model with three instruments $%
\left( Z_{1},Z_{2},Z_{3}\right) $:
\begin{equation*}
\begin{split}
Y& =\theta _{0}+\theta _{1}X+e_{1}, \\
X& =f(Z_{1},Z_{2},Z_{3})+e_{2}, \\
Z_{1}& \sim U[0,1],\qquad Z_{2}\sim U[0,1],\quad \text{ and }\quad Z_{3}\sim
\mathcal{B}(0.4),
\end{split}%
\end{equation*}%
where $(e_{1},e_{2})^{\prime }$ are non-Gaussian and correlated. Thus, $X$
in the outcome model is correlated with the error $e_{1}$. Let true $\theta
=\left( \theta _{0},\theta _{1}\right) $ equal $(1,1)$. Moreover, suppose
that the $Z_{j}$'s are relevant instruments, that is, $cov(X,Z_{j})\neq 0$
for $j\leq 3$, and
\begin{equation}
f(Z_{1},Z_{2},Z_{3})=6\left( \sqrt{0.3}Z_{1}+\sqrt{0.7}Z_{2}\right) ^{3}(1-%
\sqrt{0.3}Z_{1}-\sqrt{0.7}Z_{2})Z_{3}+Z_{1}Z_{2}(1-Z_{3}).
\end{equation}%
In practice, some instruments can be valid and some not, and the goal is to
select the valid instruments. To this end, we generate $(e_{1},e_{2},Z_{1})$
from a Gaussian copula whose covariance matrix is $\Sigma
=[1,0.7,0.7;0.7,1,0;0.7,0,1]$ 
such that the marginal distribution of $e_{1}$ is the skewed mixture of two
normal distributions $0.5\mathcal{N}(0.5,0.5^{2})+0.5\mathcal{N}(-0.5,1.118^{2})$ and the marginal distribution of $e_{2}$ is $\mathcal{N}(0,1)$. Under this setup, $Z_{1}$ is an invalid instrument. Consider the
following three models
\begin{align}
\mathcal{M}_{1}& :\mathbf{E}^{P}[(Y-\theta _{0}-\theta
_{1}X)|Z_{1},Z_{2},Z_{3}]=0, \\
\mathcal{M}_{2}& :\mathbf{E}^{P}[(Y-\theta _{0}-\theta _{1}X)|Z_{1},Z_{3}]=0,
\\
\mathcal{M}_{3}& :\mathbf{E}^{P}[(Y-\theta _{0}-\theta _{1}X)|Z_{2},Z_{3}]=0.
\end{align}%
Because $Z_{1}$ is an invalid instrument, models $\mathcal{M}_{1}$ and $%
\mathcal{M}_{2}$ are misspecified.

\noindent In $\mathcal{M}_{1}$, the basis matrix $B$ is made from the variables $%
\left( \boldsymbol{z}_{1},\boldsymbol{z}_{2},\boldsymbol{z}_{1}\odot%
\boldsymbol{z}_{2},\boldsymbol{z}_{1} \odot\boldsymbol{z}_{3},\boldsymbol{z}%
_{2}\odot\boldsymbol{z}_{3}\right)$, each using $K$ knots, concatenated
with the vector $\boldsymbol{z}_{3}$. In $\mathcal{M}_{2}$,
$B$ is made from the variables $\left( \boldsymbol{z}_{1},\boldsymbol{z}_{1} \odot\boldsymbol{z}_{3}\right)$, each using $K$ knots, concatenated with the vector $\boldsymbol{z}_{3}$. In $\mathcal{M}_{3}$, $B$ is made from the variables $\left( \boldsymbol{z}_{2},\boldsymbol{z}_{2} \odot\boldsymbol{z}_{3}\right)$, each using $K$ knots, concatenated with the vector $\boldsymbol{z}_{3}$. The number of columns in the
$B$ matrix is $ 5(K-1)+2$ for $\mathcal{M}_{1}$, and $2(K-1)+2$ for $\mathcal{M}_{2}$ and $\mathcal{M}_{3}$. The prior for $\theta _{0}$
and $\theta _{1}$ is the product of student-$t$ distributions with mean
zero, dispersion 5, and degrees of freedom equal to 2.5.
A repeated sampling experiment is conducted. The marginal likelihood of each
model is calculated in 200 repeated samples. Table \ref{tab_iv_mdd} reports
the model selection frequency for $(n=100,K =4)$, $(n=250,K=5)$, and $(n=1000,K=6)$, where
$K$ is based on $K = 2 n^{1/6}$. Note that the model
with the valid instruments, \textit{i.e.}, $\mathcal{M}_{3}$, is selected more frequently
as the number of observation gets larger, in conformity with the theory.
\begin{table}[h]
\centering
    \begin{tabular}{lccc}
    \toprule
          & $\mathcal{M}_{1}$ & $\mathcal{M}_{2}$ & $\mathcal{M}_{3}$ \\
    \midrule
    $n=100$ & 0\%     & 52\%    & 48\% \\
    $n=250$ & 0\%     & 40\%    & 60\% \\
    $n=1000$ & 0\%     & 2\%     & 98\% \\
    \bottomrule
    \end{tabular}%
\par
\vspace{0.1in} \raggedright
\caption{Model comparison: IV regression example. Each entry in the table presents
the model selection frequency in 100 repetitions; $(n=100,K =4)$, $(n=250,K=5)$, and $(n=1000,K=6)$, where $K$ is based on $K = 2 n^{1/6}$.
Each result from 10,000 MCMC draws beyond a
burn-in of 1000.}
\label{tab_iv_mdd}
\end{table}
\end{example2}

\section{Additional topics}\label{s_Additional_Topics}

\subsection{High dimensional $Z$}

We now consider the case where $Z$ lies in a high-dimensional space. If all
the elements of $Z$ are relevant, then the situation can be challenging, but there
is an interesting sub-case that is worth discussing. Suppose that the conditional expectation
depends only on a few elements of $Z$ or, in other words, most of the elements
of $Z$ are redundant. In this case, one can find the relevant elements of $Z$
by estimating and comparing models that condition on
different subsets of $Z$, where the cardinality of these subsets is say 2 or
3. The relevant elements of $Z$ correspond to the model with the largest
marginal likelihood. We refer to this procedure as sparsity-based model selection.
The next example provides an illustration.

\begin{example3}
(Sparsity-based model selection) Recall our Example 2, but assume that
one has nine additional potential $Z$'s
\begin{equation*}
Z_{j}= \frac{9}{10} Z_{1} + \frac{1}{10}\eta _{j},\quad \eta _{j}\sim \mathrm{Unif}([0,1])
\end{equation*}%
for $j=4,5,...,12$. Recall, $Z_{1}$ is an invalid instrument. Therefore, $%
Z_{j}$'s for $j=4,...,12$ are also invalid. Suppose that $Z_{3}$
affects the conditional expectation, but that one is unsure about the
remaining elements of $Z$. Suppose one believes that at most three elements of
$Z$ affect the conditional expectation (the sparsity assumption). In this
situation one can compute marginal likelihoods of the following 66 models:
\begin{equation}
 \mathbf{E}^{P}[(Y-\theta _{0}-\theta
_{1}X)|Z_{j},Z_{3}]=0 \;,  \; \; j\in \{1,2,4,5,...,12\} \;,
\end{equation}
and
\begin{equation}
\mathbf{E}^{P}[(Y-\theta _{0}-\theta
_{1}X)|Z_{j},Z_{k},Z_{3}]=0 \;, \;\; j,k\in \{1,2,4,5,...,12\} \; \text{and} \; k\neq j \;,
\end{equation}
with the correct model given by
\begin{equation}
\mathbf{E}^{P}[(Y-\theta _{0}-\theta _{1}X)|Z_{2},Z_{3}]=0
\end{equation}
Sample data (size $n=250$) is generated from the design in Example 2.
Estimation and marginal likelihood computations are based on expanded
moments from $K=3$ basis functions for each conditioning $Z_j$. A summary of
the marginal likelihood results appears in Figure \ref{fig_ex02_sparsity}, sorted by the size of the marginal likelihood. The top ranked model is the true model. As shown in the
right panel of the same figure, which reports the posterior model probabilities (under a uniform
prior on model space), the support for the true model is decisive.

\begin{figure}[h]
\begin{minipage}[c]{0.6\textwidth}
        \includegraphics[width=0.99\textwidth]{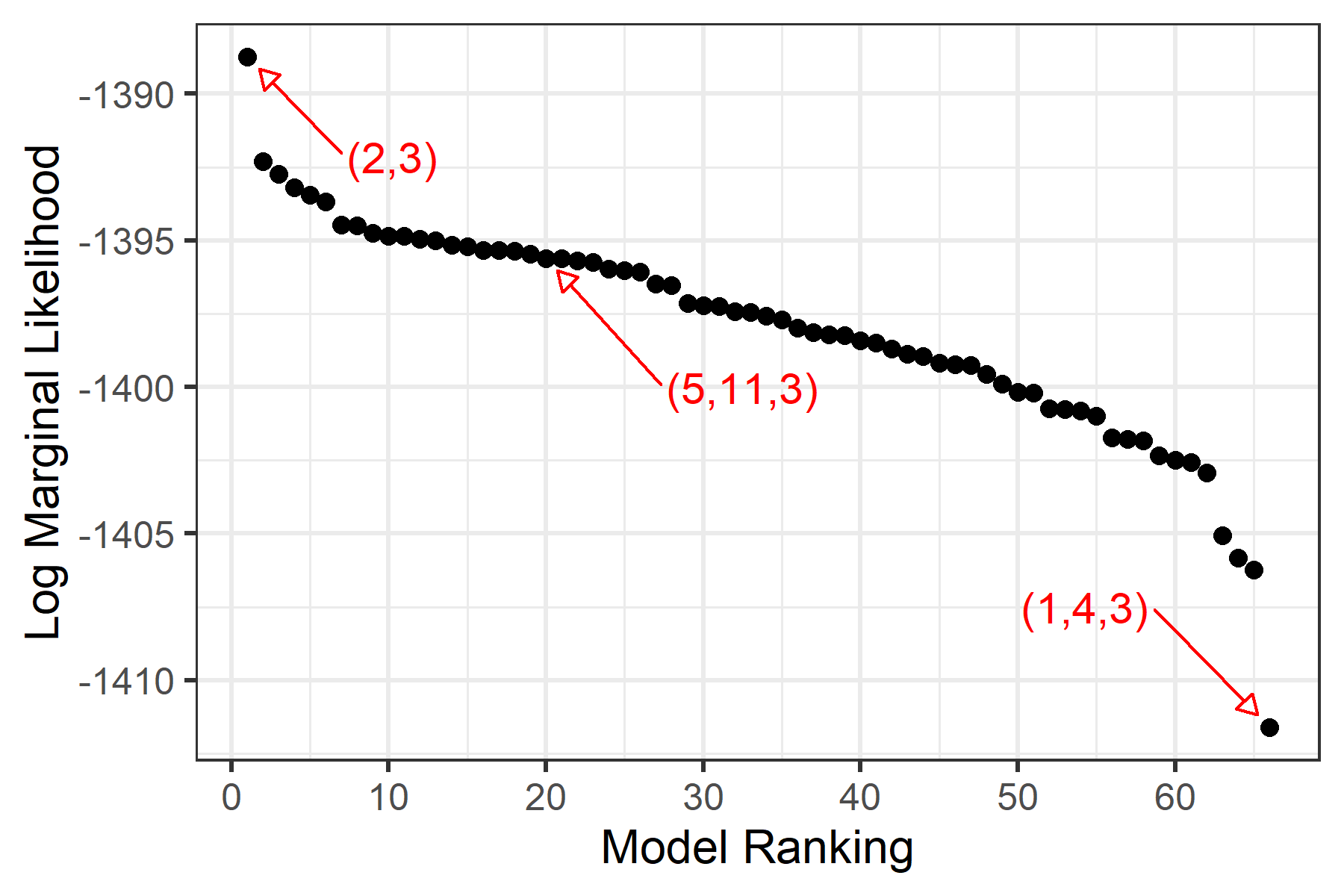}
\end{minipage}
~~
\begin{minipage}[c]{0.4\textwidth}
    \begin{tabular}{cccc}
    \toprule
    Ranking & Model & log(ML) & Prob \\
    \midrule
    1     & (2,3) & -1388.73 & 0.903 \\
    2     & (4,3) & -1392.31 & 0.025 \\
    3     & (5,3) & -1392.75 & 0.016 \\
    4     & (7,3) & -1393.20 & 0.010 \\
    5     & (5,9,3) & -1393.46 & 0.008 \\
    \bottomrule
    \end{tabular}%
\end{minipage}
\caption{Left figure presents log marginal likelihood for each of 66 models in the model space. Arrows point to the best model $(Z_{2}, Z_{3})$, top 20 model $(Z_{6}, Z_{12}, Z_{3})$, and the worst model $(Z_{1}, Z_{4}, Z_{3})$. Right table presents log marginal likelihood and posterior model probability for top 5 models. Posterior model probabilities are computed with a uniform prior on model space.}
\label{fig_ex02_sparsity}
\end{figure}

\end{example3}

\subsection{High dimensional $\protect\theta$: TaRB-MH}

It is also important to consider the estimation of conditional moment models
that contain a high-dimensional $\theta$. While the regularizing role of the
Bayesian prior is important, MCMC sampling of the posterior simulation
becomes more complicated. From our experience, the one-block M-H algorithm
that we have used above tends to be inefficient.
A straightforward alternative sampling scheme is offered by the
Tailored Randomized Block MH algorithm of \cite{ChibRamamurthy2010}. This
algorithm, which has proved useful in several similarly complex settings,
trades more computations for gains in simulation efficiency.

\begin{contexample4}
(IV regression with additional exogenous regressors). Consider the
previous IV regression model, but now with 18 additional
exogenous regressors $w$
\begin{equation*}
y_{i} = \theta_{0} + \theta_{1} x_{i} + w_{i}^{\prime }\gamma + e_{1,i}
\end{equation*}
where $w_{i} = [w_{i}^{(1)'}, w_{i}^{(2)'}, w_{i}^{(3)'}]'$, and each group $w_{i}^{(j)'}$, $j \leq 3$, are identically and independently drawn from $\mathcal{N}_6(0, \Sigma(\rho))$,
where $\Sigma(\rho)$ is a $6 \times 6$ matrix in correlation form with each off-diagonal element set
equal to 0.97. In addition, $\gamma$ is a vector of ones. In total, there are
20 unknown parameters. Other elements of the DGP are unchanged.
Suppose one has 1500 observations from this DGP, and we
estimate $[\theta_{0}, \theta_{1}, \gamma']'$ from $\mathbf{E}^{P}[(Y-\theta _{0}-\theta _{1}X)|Z_{2},Z_{3}, W]=0$ with the expanded moment conditions similar to $\mathcal{M}_{3}$ in Example \ref{example2}. The basis function matrix is formed with $Z_{2}$, $Z_{3}$, and $Z_{2} Z_{3}$, concatenated with columns in $W$. We set $K = 6$, following our recommendation, which leads to 36 expanded moment conditions. The training sample prior is based on the first 10\% of the sample, and estimation on the remaining 90\%. The prior is a product of independent student-$t$ distributions with 5 degrees of freedom, centered on the two-stage least squares (2SLS) estimate, and scale
equal to two times the 2SLS standard error.

\begin{table}[h!]
\centering
\begin{tabular}{ccccccc}
\toprule & \multicolumn{3}{c}{\underline{TaRB-MH}} & \multicolumn{3}{c}{%
\underline{One-block-MH}} \\
& Mean & SD & Ineff & Mean & SD & Ineff \\
\midrule
    $\theta_{0}$ & 1.03  & 0.03  & 4.35  & 1.03  & 0.03  & 14.34 \\
    $\theta_{1}$ & 0.91  & 0.11  & 7.57  & 0.91  & 0.10  & 47.38 \\
    $\gamma_{1}$ & 0.99  & 0.03  & 6.92  & 0.99  & 0.03  & 13.83 \\
    $\gamma_{2}$ & 0.93  & 0.15  & 6.06  & 0.92  & 0.15  & 15.09 \\
    $\gamma_{3}$ & 0.96  & 0.15  & 2.75  & 0.97  & 0.15  & 12.41 \\
    $\gamma_{4}$ & 1.14  & 0.17  & 2.21  & 1.15  & 0.17  & 12.54 \\
    $\gamma_{5}$ & 1.32  & 0.16  & 1.97  & 1.33  & 0.16  & 14.26 \\
    $\gamma_{6}$ & 1.02  & 0.16  & 1.57  & 1.02  & 0.15  & 12.79 \\
    $\gamma_{7}$ & 0.99  & 0.03  & 7.05  & 0.99  & 0.02  & 10.58 \\
    $\gamma_{8}$ & 1.04  & 0.14  & 5.42  & 1.04  & 0.14  & 14.60 \\
    $\gamma_{9}$ & 1.18  & 0.15  & 2.98  & 1.18  & 0.14  & 11.73 \\
    $\gamma_{10}$ & 1.22  & 0.16  & 2.37  & 1.22  & 0.15  & 14.66 \\
    $\gamma_{11}$ & 0.94  & 0.15  & 1.81  & 0.94  & 0.15  & 11.56 \\
    $\gamma_{12}$ & 0.64  & 0.16  & 1.69  & 0.64  & 0.16  & 24.87 \\
    $\gamma_{13}$ & 1.01  & 0.03  & 7.00  & 1.01  & 0.03  & 12.89 \\
    $\gamma_{14}$ & 0.86  & 0.15  & 5.58  & 0.86  & 0.15  & 14.66 \\
    $\gamma_{15}$ & 0.84  & 0.15  & 3.34  & 0.84  & 0.15  & 14.86 \\
    $\gamma_{16}$ & 1.04  & 0.16  & 2.24  & 1.05  & 0.15  & 12.84 \\
    $\gamma_{17}$ & 1.18  & 0.17  & 2.11  & 1.18  & 0.17  & 20.25 \\
    $\gamma_{18}$ & 0.86  & 0.16  & 1.44  & 0.86  & 0.16  & 16.31 \\
\bottomrule
\end{tabular}%
\caption{Posterior summary of IV regression example with additional
covariates ($n=1500$). The true value of all parameters ($\protect\theta$'s
and $\protect\gamma$'s) are set to one. The summaries are based on 50,000
MCMC draws beyond a burn-in of 10,000 for the one-block-MH sampler and 3,000
draws beyond a burn-in of 1,000 for the TaRB-MH. The M-H acceptance rate is around 37\% for the one-block-MH and 87\% for TaRB-MH. ``Ineff'' is the inefficiency factor.}
\label{tab:iv_tarbmh}
\end{table}
\end{contexample4}

\noindent The results appear in Table \ref{tab:iv_tarbmh}. In implementing
the TaRB-MH sampling scheme, the probability of starting a new block is set to 0.3,
so that the each block, within each MCMC iteration, contains 6 parameters on average.
For comparison, results from the single-block sampling scheme (on the
same conditional moments) are also included. It is evident that the
two MCMC samplers produce identical posterior moments, but that the TaRB-MH
sampler dominates the one-block MH sampler in terms of simulation efficiency
as measured by the inefficiency factor (the ratio of the numerical variance
of the mean to the variance of the mean assuming independent draws). An
inefficiency factor close to 1 indicates that the MCMC draws
are essentially independent. Therefore, armed with the TaRB-MH sampler,
computational efficiency is retained, even in higher-dimensional $\theta$
problems.

\section{Applications}\label{s_Applications}

\subsection{Asset pricing}

A key question in finance concerns the makeup of the pricing kernel, or the
stochastic discount factor (SDF). Factors in the SDF are the risk factors
that explain the cross-section of expected equity returns and, for this reason,
establishing the identity of these risk factors has been a long-standing quest in
finance.

Following notation from \cite{ChibZeng2020}, write the SDF $M_{t}$ at time
(month) $t$ as
\begin{equation}
M_{t}=1-b^{\prime }(x_{t}-\mu _{x})
\end{equation}%
where $x_{t}$ is a $(k_{x}\times 1)$ vector of risk factors (empirically these are the excess returns
on portfolios of stocks), and $b$ is the unknown
risk-factor premia and $\mu _{x}=$ $\mathbf{E}^{P}(x_{t})$.\textbf{\ }The
parameters $(b,\mu _{x})$ are unknown. Suppose that there are
other factors (excess returns on other portfolios) that are collected in a
$(k_{w}\times 1)$ vector $w_{t}$. Let $f_{t} \coloneqq(x_{t}^{\prime },w_{t}^{\prime})^{\prime }$ be a $(k_{f}\times 1)$-vector, where $k_{f}=k_{x}+k_{w}$. If $x_{t}$ are risk factors, then finance theory dictates that the restriction $\mathbf{E}^{P}(M_{t}f_{t})=0$ holds. Given a sample of observations $%
\{f_{t}\}_{t=1}^{n}$, one can estimate $(b,\mu _{x})$ based on the following
moment conditions
\begin{equation*}
\mathbf{E}^{P}[(1-b^{\prime }(x_{t}-\mu _{x}))f_{t}]=0,\qquad \mathbf{E}^{P}%
\left[ x_{t}-\mu _{x}|f_{t-1}\right] =0,
\end{equation*}
where the second conditional moment restriction identifies $\mu _{x}$.

As an example, consider the data at
\url{http://apps.olin.wustl.edu/faculty/chib/rpackages/czfactor/czfactor.pdf}
on monthly excess returns (Jan 1974 -- Dec 2018) on $k_{f}=12$ potential
risk-factors. Thus, in this situation, there are 12 conditioning variables,
an illustration of a modestly high-dimensional $Z$. Let $x_{t}$ be the
excess return on the market portfolio (denoted Mkt in the data).

Now construct the expanded moment conditions as
\begin{equation}
\mathbf{E}^{P}\left[ (x_{t}-\mu _{x})\otimes \lbrack q^{K}(f_{1,t-1}),%
\widetilde{q}^{K}(f_{2,t-1}),\ldots,\widetilde{q}^{K}(f_{12,t-1})]\right] =0,
\end{equation}%
where $q^{K}(f_{1,t-1})$ consist of $K=3$ basis functions, and $\tilde{q}%
^{K}(f_{j,t-1})$ $(j\geq 2)$ each consist of 2 basis functions derived from $%
q^{K}(f_{j,t-1})$ by subtracting the second and third columns from the first
and then dropping the first. Along with these 25 expanded moment conditions,
the pricing conditions supply an additional twelve, for a total of 37 moment
conditions.

For the prior, one can employ the training sample approach. From the first
80 observations (the training sample) the hyperparameters of the independent
student-t distribution of $(b,\mu _{x})$ with 2.5 degrees of freedom are set
as follows. The center of the prior density is set to the Generalized Method
of Moments (GMM) estimate, and the scale to two times the GMM\ standard
error. This black-box prior is particularly convenient if the analysis has
to be repeated for different possible variables in the SDF. The remaining 459
observations are used to construct a joint posterior distribution of $b$ and
$\mu_{x}$.

\begin{figure}[h!]
\centering
\includegraphics[width=4in]{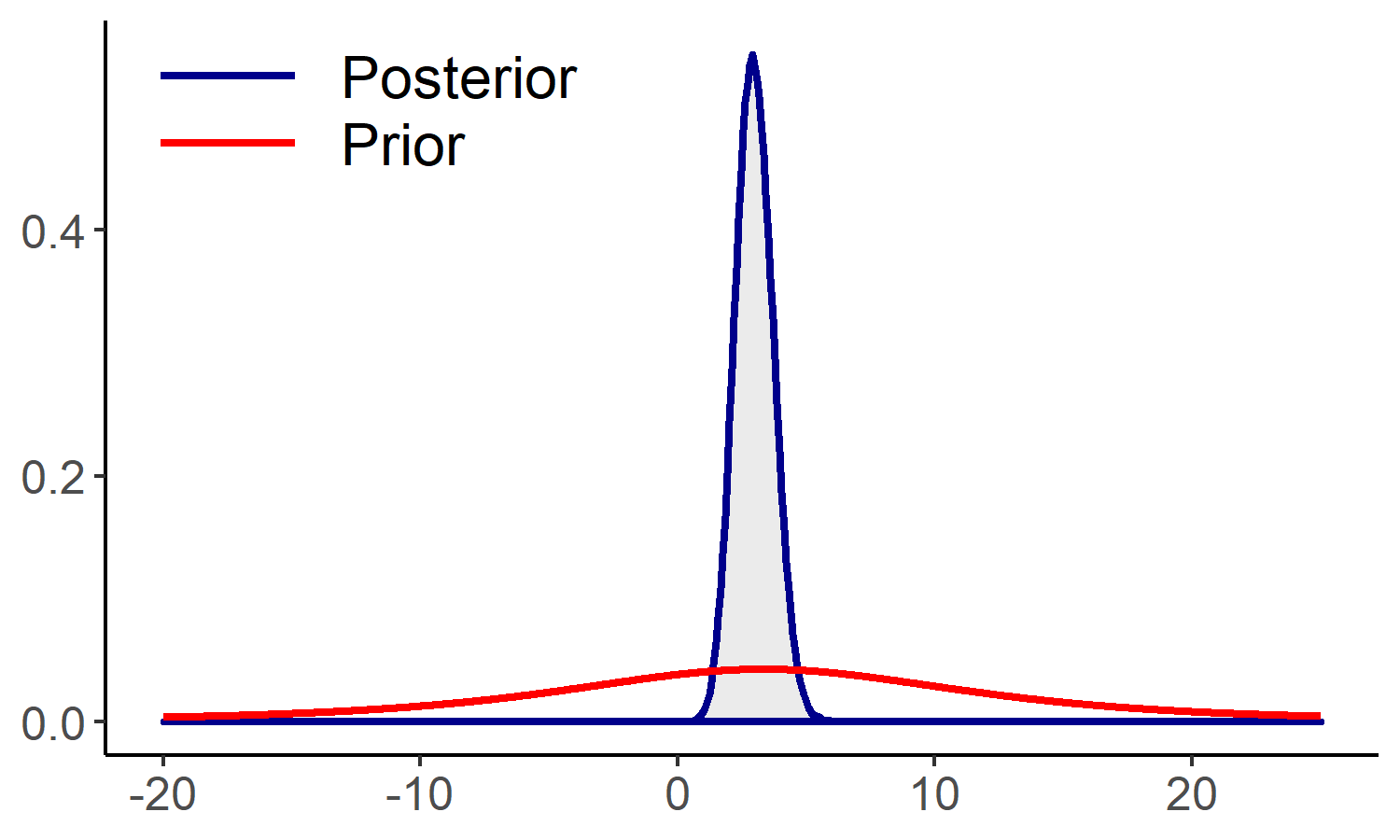}
\caption{Asset pricing data: Prior and posterior density of $b$.
Prior is Student-$t$ density based on the first 80
observation. Posterior density is based on the remaining 459 observations.
50,000 MCMC draws after 1,000 burn-in.}
\label{fig_finance_example}
\end{figure}

\begin{table}[h!]
\centering
\begin{tabular}{lccccc}
\toprule & Mean & SD & Median & 5\% & 95\% \\
\midrule
    $b$   & 2.981 & 0.730 & 2.955 & 1.818 & 4.211 \\
    $\mu_{x}$ & 0.006 & 0.001 & 0.006 & 0.004 & 0.008 \\
\bottomrule &  &  &  &  &
\end{tabular}%
\caption{Asset pricing data: Summary of the posterior distribution
based on 50,000 MCMC draws after 1,000 burn-in.}
\label{tab:risk-factor1}
\end{table}
The posterior summary of $(b,\mu_{x})$ from 50,000 MCMC draws is given in
Table \ref{tab:risk-factor1},
and the prior and posterior density of $b$ are
presented in Figure \ref{fig_finance_example}.
This summary, specifically the lower and upper limits of
the marginal posterior of $b$ confirm that $b$ is non-zero,
and hence that the Mkt variable is a risk-factor.

\subsection{ATE under conditional ignorability}

\label{s_5}

For another important application of the methods in this paper, consider the
problem of estimating the average treatment effect (ATE) under the
assumption of conditional ignorability. In the frequentist literature, the
ATE is commonly estimated by propensity score methods and, on the Bayesian
side, from models of the potential outcomes. These models generally have a
nonparametric mean function, but parametric noise.
By adopting the conditional moment perspective, however,
one can evade the burden of distributional assumptions.

Data is from the 1997 Child Development Supplement to the Panel Study of
Income Dynamics \citet[Section 5.8.2]{GuoFraser15}, where the ATE is
calculated by the propensity score. The research question is the effect of
childhood welfare dependency on academic achievement. The latter, the
dependent variable $y$, is measured by the child's score on the
\textquotedblleft letter-word identification\textquotedblright\ section of
the Woodcock-Johnson Revised Tests of Achievement. The treatment variable $x$
equals one if the child received AFDC (Aid to Families with Dependent
Children) benefits at any time from birth to 1997 (the survey year) and
equals zero if the child never received benefits during that period. It is
assumed that the potential outcomes are independent of $x$, conditioned on $%
z_{1}$,$z_{2},...,z_{6}$ (the assumption of conditional ignorability), where

\begin{itemize}
\item $z_{1}$: mratio97, the ratio of family income to the poverty line in
1997

\item $z_{2}$: pcged97, the caregiver's years of schooling

\item $z_{3}$: pcg\_adc, the number of years in which the caregiver received
AFDC in her childhood

\item $z_{4}$: age97, the child's age in 1997

\item $z_{5}$: race, one for African-American children and zero for other

\item $z_{6}$: male, one if the child is male and zero if female.
\end{itemize}

Two observations from the sample are dropped. These have values of mratio97
larger than $9$ standard deviation from the mean of mratio97. Apart from
mratio97 and age97, the other variables are categorical. There are $n_{0}=727 $
control subjects and $n_{1}=274$ treated subjects. The ATE is
expected to be negative, reflecting the hypothesis that
welfare dependency has an adverse effect on academic achievement.

To answer the research question, suppose that the potential outcomes for the
controls satisfy the conditional moments
\begin{equation*}
\mathbf{E}^{P}((y_{i0}-\beta _{00}-h_{01}(z_{1})-\beta _{02}z_{2}-\beta
_{03}z_{3}-h_{04}(z_{4})-\beta _{05}z_{5}-\beta _{06}z_{6})|z_{i})=0
\end{equation*}%
and those for the treated satisfy the conditional moments
\begin{equation*}
\mathbf{E}^{P}((y_{i1}-\beta _{10}-h_{11}(z_{1})-\beta _{12}z_{2}-\beta
_{13}z_{3}-h_{14}(z_{4})-\beta _{15}z_{5}-\beta _{16}z_{6})|z_{i})=0
\end{equation*}%
where $\{h_{01},h_{04},h_{11},h_{14}\}$ are four non-parametric functions.
These are each modeled by natural cubic splines with 5 knots. Thus, the
parameters $\theta _{j}$ of the $j$th potential outcome model consist of $%
(\beta _{j0},\beta _{j2},\beta _{j3},\beta _{j5},\beta _{j6})$ plus the
eight spline coefficients. Special cases of this model, mentioned below, are
considered and evaluated by marginal likelihoods. For example, models in
which the $h$ functions are linear are of interest.

The expanded moments are constructed as follows. The basis matrix has cubic
spline basis functions for $(z_{1},z_{4})$, each with 5 knots, concatenated
with $(z_{2},z_{3},z_{5},z_{6})$ (as is) because the latter variables are
all essentially categorical. In total, this produces 13 expanded
unconditional moments for the estimation of the $y_{0}$ and $y_{1}$ models.
The prior distribution on the parameters is a product of student-t
distributions with 2.5 degrees of freedom with mean of the intercept
equal to the mean of the first 50 observations, the mean of the slopes
equal to 0, and dispersion equal to 5.

Four models are estimated and evaluated. In the baseline model, the $h$
functions are linear. In the second model, only the effect of $z_{1}$ is
nonparametric. In the third model, only the effect of $z_{4}$ is assumed to
be nonparametric and, finally, in the fourth model, both $z_{1}$ and $z_{4}$
are nonparametric. The results given in Table~\ref{tab:academicmod} show
that the model best supported by these data is the third.

\begin{table}[h!]
\centering
\par
\begin{tabular}{lcc}
\toprule & Non-treated & Treated \\
\midrule Linear & -4823.76 & -1555.55 \\
\midrule $z_{1}$ nonparametric & -4829.23 & -1556.23 \\
\midrule $z_{4}$ nonparametric & \textbf{-4813.86} & \textbf{-1555.78} \\
\midrule $z_{1}$ and $z_{4}$ nonparametric & -4818.98 & -1556.56 \\
\bottomrule &  &
\end{tabular}%
\caption{Academic achievement data: Marginal likelihoods of 4 competing
causal models, based on 20,000 MCMC draws beyond a burn-in 1000.}
\label{tab:academicmod}
\end{table}

\begin{figure}[h!]
\begin{center}
\includegraphics[width=4in]{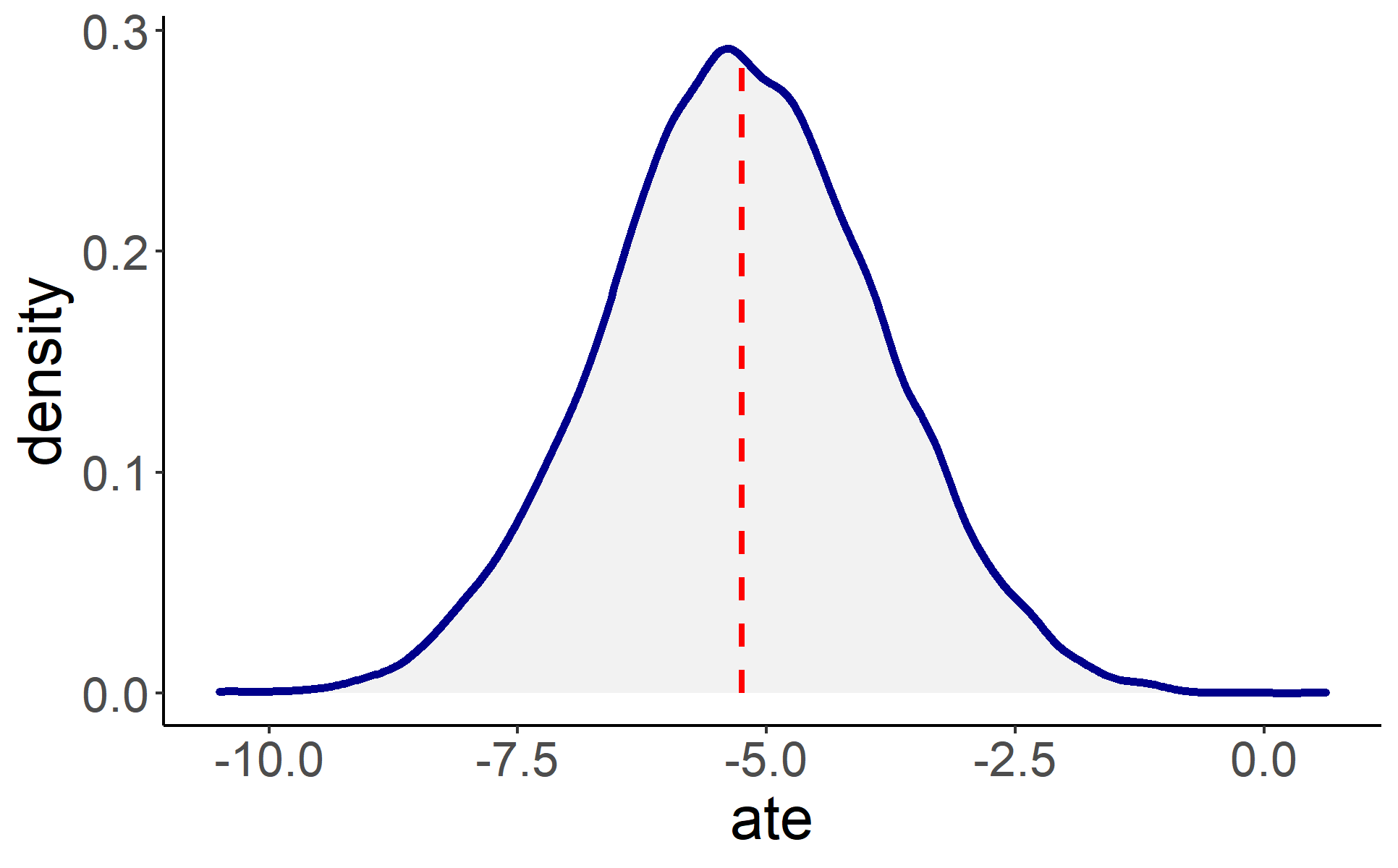}
\par
\centering
\end{center}
\caption{Academic achievement data: Posterior density of the ATE
from model in which the effect of $z_4$ is non-parametric. Red
dashed line is at -5.251, the posterior mean. Posterior density
based on 20,000 posterior draws after 1,000 burn-in.}
\label{fig:GuoFraserAte}
\end{figure}

Consider now posterior inference on the ATE. By definition,
the sample version of the ATE is
\begin{equation*}
\text{ATE}=\frac{1}{n}\sum_{i=}^{n}\left( \mathbf{E}^{P}\left(
y_{i1}|z_{i},\theta _{1}\right) -\mathbf{E}^{P}\left( y_{i0}|z_{i},\theta
_{0}\right) \right) \;,
\end{equation*}
where, in the model selected by the preceding comparison,
\begin{equation*}
\mathbf{E}^{P}\left( y_{ij}|z_{i},\theta _{j}\right) =\beta_{j0}+
\beta _{j1}z_{1}+\beta _{j2}z_{2}+\beta _{j3}z_{3}+ h_{j4}(z_{4}) + \beta
_{j5}z_{5}+\beta _{j6}z_{6} \;.
\end{equation*}%
Clearly, if we evaluate the latter expression
at each posterior draw of $\left( \theta_{0},\theta_{1}\right)$,
we produce a sample of the ATE from its posterior distribution.
We summarize this sample in Table~\ref{table:GuoFraserSATE} and
Figure~\ref{fig:GuoFraserAte}.
\begin{table}[h!]
\centering
\begin{tabular}{lrrrrr}
\hline
& \multicolumn{1}{c}{Mean} & \multicolumn{1}{c}{SD} & \multicolumn{1}{c}{
Median} & \multicolumn{1}{c}{Lower} & \multicolumn{1}{c}{Upper} \\
\hline\hline
Propensity Score Matching & -5.682 & 1.976 &  & -9.496 & -1.502 \\
Bayesian ATE & -5.257 & 1.393 & -5.267 & -7.983 & -2.545 \\ \hline
\end{tabular}%
\caption{Academic achievement data: Summary of the posterior distribution of
the ATE from model in which the effect of $z_4$ is non-parametric.}
\label{table:GuoFraserSATE}
\end{table}
One can see that the ATE posterior point estimate is
similar in size to the propensity score estimate, but
the posterior standard deviation is smaller, leading to a
less dispersed interval estimate.
As a takeaway, it is striking that the Bayesian analysis of
this important problem can be prosecuted under such
minimal assumptions.

\section{Conclusion}

\label{s_6}

In this paper we have developed perhaps the first Bayesian framework for
analyzing an important and broad class of semiparametric models in which the
distribution of the outcomes is defined only up to a set of conditional
moments, some of which may be misspecified. We have derived BvM theorems for
the behavior of the posterior distribution under both correct and incorrect
specification of the conditional moments, and developed the theory for
comparing different conditional moment models through a comparison of model
marginal likelihoods. In addition, we have discussed settings with a
high-dimensional $Z$ and $\theta$, the former addressed by a sparsity-based
model search procedure, and the latter by the TaRB-MH MCMC algorithm
for efficient posterior sampling.

Our theory and various examples, taken together, show that the developments in this
paper make possible, for the first time, the formal (and practical) Bayesian analysis
of a new, large class of problems that were hitherto difficult, or not possible, to
tackle from the Bayesian viewpoint. This research we believe
should have numerous positive ramifications for the growth and practice and teaching of
Bayesian statistics.

\subsection*{Supplementary Materials}

The supplementary materials in the online Appendix contain the technical proofs of the results in the paper, additional applications, and R code, packages and data for replicating the examples. This replication content can also be accessed from
\url{https://apps.olin.wustl.edu/faculty/chib/chibshinsimoni2021replication/}.
%
\setlength{\baselineskip}{19.0pt}

\setlength{\baselineskip}{21.0pt}

\setlength{\bibsep}{1.5pt}

\bibliographystyle{plain}
\bibliography{AnnaBib}


\end{document}